\documentclass[a4paper,12pt]{amsart} 

\usepackage{epsfig}
\usepackage{graphicx}

\usepackage[latin1]{inputenc}
\usepackage[T1]{fontenc}
\usepackage{indentfirst}
\usepackage{amssymb}
\usepackage{eufrak}
\usepackage{amsmath}
\usepackage{amsfonts}
\usepackage{amsthm}
\usepackage{mathrsfs}
\usepackage[bookmarks]{hyperref}
\usepackage{xypic}
\newcommand{\R}{\mathop{\rm Re}\nolimits}

\newcommand{\Log}{\mathop{\rm Log}\nolimits}
\newcommand{\val}{\mathop{\rm val}\nolimits}
\newcommand{\Val}{\mathop{\rm Val}\nolimits}
\newcommand{\ord}{\mathop{\rm ord}\nolimits}
\newcommand{\Arg}{\mathop{\rm Arg}\nolimits}
\newcommand{\Vol}{\mathop{\rm Vol}\nolimits}
\newcommand{\Ima}{\mathop{\rm Im}\nolimits}

\newcommand{\Area}{\mathop{\rm Area}\nolimits}
\newcommand{\supp}{\mathop{\rm supp}\nolimits}
\newcommand{\Critv}{\mathop{\rm Critv}\nolimits}
\newcommand{\Verte}{\mathop{\rm Vert}\nolimits}

\input epsf

\def \square{\smallskip \hfill \vrule width 5 pt height 7
pt depth - 2 pt \smallskip }

\newenvironment{prooof}
{\noindent {{\it Proof} \;}}{\hspace*{\fill}\square\vskip 8pt}

\theoremstyle{plain}
\newtheorem{Lem}[subsection]{Lemma}
\newtheorem{The}[subsection]{Theorem}
\newtheorem{Cor}[subsection]{Corollary}
\newtheorem{Pro}[subsection]{Proposition}
\theoremstyle{definition}
\newtheorem{Rem}[subsection]{Remark}
\newtheorem{Rems}[subsection]{Remarks}

\newtheorem{Def}[subsection]{Definition}
\newtheorem{Exe}[subsection]{Example}

\begin{document}

\title
{Geometric and Combinatorial Structure of Hypersurface Coamoebas}


\author{Mounir Nisse}
\address{Université Pierre et Marie Curie - Paris 6, IMJ (UMR 7586),
Labo: Analyse Algébrique, Office: 7C14, 175, rue du Chevaleret,\\ 75013
Paris, France}
\email{nisse@math.jussieu.fr}

\maketitle

\begin{abstract} Let $V$  be a complex algebraic hypersurface defined by a polynomial $f$ with Newton polytope $\Delta$. It is well known that the spine of its amoeba has a structure of a tropical hypersurface. We prove in this paper that there exists a complex tropical hypersurface $V_{\infty ,\, f}$ such that its coamoeba is homeomorphic to the closure in the real torus of the coamoeba of $V$. Moreover, the coamoeba of $V_{\infty ,\, f}$  contains an arrangement of $(n-1)$-torus depending only on the geometry of $\Delta$ and the  coefficients of $f$. In addition, we can consider this arrangement, as a weighted codual hyperplanes arrangement in the universal covering of the real torus, and the balancing condition (the analogous to that of tropical hypersurfaces) is satisfied. This codual hyperplanes arrangement is called the {\em shell} of the complex coamoeba (the cousin of the spine of the complex amoeba).
Using this combinatorial coamoebas structure, we show that the amoebas of complex algebraic hypersurfaces defined by maximally sparse polynomials are solid. More precisely, we characterize the image of the order map defined by Forsberg, Passare, and Tsikh.  
\end{abstract}

\setcounter{tocdepth}{1} \tableofcontents

\section{Introduction}

Amoeba and coamoeba are a very  fascinating notions in mathematics where the first  terminology has been introduced by I. M. Gelfand, M M. Kapranov and A. V. Zelevinsky in their book (see \cite{GKZ-94}) in 1994,
and the second one by M. Passare and A. Tsikh in 2001. Amoebas (resp. coamoebas) have their spines, contours and tentacles (resp. spines, contours and extra-pieces), and they have many applications in real algebraic geometry , complex analysis, mirror symmetry and in several other areas (see \cite{M1-02}, \cite{M2-04}, \cite{M3-00}, \cite{FPT-00}, \cite{PR1-04}, \cite{R-01}, \cite{RST-05},  \cite{S-02}, and \cite{SS-04}). Amoebas and coamoebas are naturally linked to the geometry of Newton polytopes, which can be seen in particular with the Viro patchworking principle (i.e., tropical localization) based on the combinatorics of subdivisions of convex lattice polytopes.
The purpose of this paper is to describe the relations and the similarities which exist between amoebas and coamoebas of algebraic complex hypersurfaces.
Let $V\subset (\mathbb{C}^*)^n$ be an algebraic  complex hypersurface defined by a polynomial $f$ with Newton polytope $\Delta$. The amoeba $\mathscr{A}$ of an algebraic set $V=\{ f(z)=0 \}$ in the algebraic torus $(\mathbb{C}^*)^n$ is defined as its image under the mapping $\Log : (z_1,\ldots , z_n)\mapsto (\log \Vert z_1\Vert ,\ldots , \log \Vert z_n\Vert )$. The amoeba's complement has a finite number of convex connected components, corresponding to domains of convergence of the Laurent series expansions of the rational function $\frac{1}{f}$.
We know that the spine $\Gamma$ of the amoeba $\mathscr{A}$ has a structure of a tropical hypersurface in $\mathbb{R}^n$ (proved by M. Passare and H. Rullg\aa rd in 2000 \cite{PR1-04}, and independently by G. Mikhalkin in 2000). In addition the spine of the amoeba is dual to some coherent (i.e. convex) subdivision $\tau$ of the integer convex polytope $\Delta$. It is shown by M. Forsberg, M. Passare and A Tsikh that the set of vertices of $\tau$ is in bijection with the set of complement components of $\mathscr{A}$  in $\mathbb{R}^n$ \cite{FPT-00}. The coamoeba $co\mathscr{A}$ of an algebraic set $V=\{ f(z)=0\}$ in $(\mathbb{C}^*)^n$ is defined as its image under the argument mapping $\Arg : (z_1,\ldots ,z_n)\mapsto (e^{i\arg (z_1)},\ldots ,e^{i\arg (z_n)})$. It is shown in \cite{N1-07} that the complement components of the closure in the flat torus of the coamoeba of a complex algebraic hypersurface defined by a polynomial $f$ with Newton polytope $\Delta$ are convex and their number cannot exceed $n!\Vol (\Delta )$.

\vspace{0.1cm}

The study of the geometry and topology of the amoeba of an algebraic sub-variety in the complex 
algebraic torus is strongly linked to the coamoeba of that sub-variety. I am convinced that these two 
objects are not only similar, but also, complementary. Also, their combinatorial structure 
is in duality with an object defined by the defining ideal of  this sub-variety. For example,
such complementarity is proved in codimension one. Moreover, the combinatorial structure of
the Newton polytope in this case is in duality with that of the amoeba and coamoeba. 
More precisely, these two objects are related combinatorially to the same object which is some convex (coherent)
subdivision $\tau$ of the Newton polytope. The amoeba has a spine $\Gamma$ (which can be seen as its average) equipped
with a tropical structure dual to the subdivision $\tau$, and the coamoeba has a {\em shell} which is a codual 
hyperplanes arrangement $\mathscr{H}$
(codual to the edges of $\tau$) equipped with a structure of a co-tropical structure. I mean that $\mathscr{H}$ is
a weighted polyhedral complex satisfying to the balancing condition as the spine is. So, this arrangement
 has an algebraic, topological, and combinatorial structure strongly linked to the geometry and the topology 
 of the complex hypersurface itself. In this paper, we prove the existence  of the arrangement $\mathscr{H}$ by giving an explicit
  function $\nu_{(2,\, f)}$ defined
 on a subset of the support of the polynomial defining the complex hypersurface union the image of the order map,
  which take its image in $S^1$. This function gives a complete description of the arrangement $\mathscr{H}$.
  A combinatorial structure of coamoebas of complex plane curves is given by physicists, and described by bipartite graphs
  called {\em brane tiling} see \cite{FHKV-05} and \cite{HV-07}, but their description has no natural 
  generalizations in higher dimension, maybe in the three
  dimensional case for a very special Newton polytopes. Our combinatorial description is naturally and strongly related to the
  Newton polytope and the coefficients of the defining polynomial of the complex hypersurface.

Sturmfels's tropical model, I mean the image under the valuation map of an algebraic 
sub-variety in the algebraic torus over the field of Puiseux series $\mathbb{K}$ is the standard model
 in tropical geometry, which plays the central role in tropical geometry.
So, to study the amoeba (resp. coamoeba) of an algebraic sub-variety in the complex 
algebraic torus of codimension greater than one, we must firstly study those amoebas (resp. coamoebas)
which converge in the Gromov-Hausdorff metric, to Sturmfels's tropical (resp. co-tropical) model.
 Here, Sturmfels's co-tropical model means the coamoeba of some lifting in the complex  algebraic torus,
 of a tropical model (see section 3 for more details).
 It may be noted that in codimension one, any tropical variety is of Sturmfels's model, but in codimension 
greater than one it can not be globally but only locally modeled. In the forthcoming papers with F. Sottile 
\cite{NS1-09},
we study the combinatorial structure of coamoebas of some complex tropical complete intersections, and also
the coamoebas of some complex tropical Grassmannian.

\vspace{0.2cm}

The purpose of this paper is to prove that the coamoebas of a complex algebraic hypersurfaces have a similar combinatorial object as the spine of their amoebas.  I also, give a necessary condition  of  combinatorial nature, which must satisfies  an index in the image of the order map but not in the support of the polynomial.

\vspace{0.2cm}

\begin{The}  Let $V$ be a complex algebraic hypersurface defined by a polynomial $f$ with Newton polytope $\Delta$. Let us denote by $\tau_f$ the subdivision of $\Delta$ dual to the spine of the amoeba of $V$.
 Then there exists a complex tropical hypersurface $V_{\infty ,\, f}$ satisfying the following:
\begin{itemize}
\item[(i)]\,  The  closure of the coamoebas of $V_{\infty ,\, f}$ and $V$ in the real torus $(S^1)^n$ are homeomorphic;
\item[(ii)]\, The lifting of the coamoeba of $V_{\infty ,\, f}$ in the universal covering of the torus $(S^1)^n$ contains an arrangement $\mathscr{H}$ of codual  hyperplanes to the set of edges of $\tau_f$ which determine completely the topology of the complex coamoeba of $V$.
\end{itemize}
\end{The}

\vspace{0.1cm}

Let us begin with a brief description of our ideas without technical details. Let $A$  be the support of the polynomial $f$. The main ingredients in the construction are a special  deformation of the standard complex structure on $(\mathbb{C}^*)^n$, the generalized Passare-Rullg\aa rd function $\nu$ (see the definition in section 2), Viro's  tropical localization, and  Kapranov's and Sturmfels's Theorem   \cite{K-00}, \cite{S-02}, and \cite{SS-04}. 
So, 
 using the generalized Passare-Rullg\aa rd function $\nu$, we construct a family of polynomials $f_t$ with $0<t\leq \frac{1}{e}$ (i.e., a deformation) such that $f_{\frac{1}{e}} = f$, and we consider the family of the $J_t$-holomorphic hypersurfaces $H_t(\{ f_t(z) = 0 \} )$ where $H_t$ is a self-diffeomorphism of $(\mathbb{C}^*)^n$. When $t$ tends to zero, we obtain a complex tropical hypersurface $V_{\infty ,\, f}$, such that its coamoeba is a retract by deformation of the coamoeba of $V$.  In addition, using the subdivision $\tau$ of $\Delta$ dual to the spine of the amoeba of $V$ and Sturmfels's theorem \cite{S-02} and \cite{SS-04}, that is the several ways to view a tropical variety, or Kapranov's theorem \cite{K-00},  we have an algorithm giving an explicit description of the coamoeba of $V_{\infty ,\, f}$.
 In other word, the results are obtained by deformation of the complex structure on the hypersurface to a degenerate structure called complex tropical structure  which is a  piecewise-linear polyhedral complex in $\mathbb{R}^n$ supplied with some lifting to $(\mathbb{C}^*)^n$ (see G. Mikhalkin \cite{M1-02} and \cite{M2-04}).

\vspace{0.2cm}

\begin{The} Let $\beta\in \Ima (\ord_f)\setminus \supp (f) $. Then we have the following:
\begin{itemize}
\item[(i)]\, There exists an effective sub-arrangement $\mathscr{H}'$ in $\mathscr{H}_f$ of codual hyperplanes to the adjacent edges to $\beta$
such that any hyperplane $H\in \mathscr{H}'$ is of weight at least two;
\item[(ii)]\, The coamoeba $co\mathscr{A}$ of the hypersurface $V$ with defining polynomial $f$ contains some extra-pieces. In particular, if $n=2$, this imply that the real part $\mathbb{R}V$ of $V$ is not a Harnack curve.
\end{itemize}
\end{The}

\vspace{0.1cm}

\begin{The} The amoeba of a complex algebraic hypersurface  defined by  maximally sparse polynomial is solid.
\end{The}

\vspace{0.1cm}

\begin{The} Let $V$ be a complex  algebraic  plane curve defined by a polynomial $f$ with Newton polytope $\Delta$, such that its real part $\mathbb{R}V$ is a Harnack curve. Then $f$ is dense i.e., $\supp (f) = \Delta\cap\mathbb{Z}^2$.
\end{The}

\vspace{0.1cm}

\noindent Recall that it is proved by M. Passare, T. Sadykov, and A. Tsikh  \cite{PST-05} that the amoeba of 
any $A$-discriminantal hypersurface is solid.

\vspace{0.2cm}
The remainder of this paper is organized as follows. In  Section 2, we review some  properties of the amoebas of  complex hypersurfaces proved by  M. Forsberg, M. Passare  and A. Tsikh in \cite{FPT-00}, M Passare and Rullg\aa rd in \cite{PR1-04}, and G. Mikhalkin in \cite{M1-02} and \cite{M2-04}. Also, it reviews some theorem structure of non-Archimedean amoebas proved by M. Kapranov in \cite{K-00}, B. Sturmfels in \cite{S-02}, and D. Speyer and B. Sturmfels in \cite{SS-04}. It also reviews some properties of the coamoebas of complex hypersurfaces proved in \cite{N1-07} and \cite{N2-07}. In Section 3, we review some properties of complex tropical hypersurfaces. In Section 4, we define the codual hyperplanes arrangement in the universal covering of the real torus with other definitions, properties, and examples. In Section 5, we give the proof of the main result i.e., Theorem 1.1 which gives a geometric and combinatorial structure of complex hypersurface coamoebas. In Section 6, we give a combinatorial caracterization of the lattice points in the image of the order map $\ord$, and then we give a second proof that the amoeba of maximally sparse polynomial is solid. We will prove also, that if a real polynomial $f$  in two variables with Newton polygon $\Delta$ is the defining polynomial of a Harnack curve, then it is dense i.e., $\supp (f) = \Delta\cap \mathbb{Z}^2$.

\section{Preliminaries}

In this paper we will consider  algebraic hypersurfaces $V$ in the
complex torus $(\mathbb{C}^*)^n$, where $\mathbb{C}^*
=\mathbb{C}\setminus \{ 0\}$ and $n\geq 1$ an integer. This means that $V$ is the zero locus of a polynomial:
$$
f(z) =\sum_{\alpha\in \supp (f)} a_{\alpha}z^{\alpha},\,\,\,\,\quad\quad\quad
z^{\alpha}=z_1^{\alpha_1}z_2^{\alpha_2}\ldots z_n^{\alpha_n}
\quad\quad\quad\quad\quad\quad\quad\quad  (1)
$$
where each $a_{\alpha}$ is a non-zero complex number and $\supp (f)$ is a
finite subset of $\mathbb{Z}^n$, called the support of the polynomial
$f$, with convex hull, in $\mathbb{R}^n$, the Newton polytope
$\Delta_f$ of $f$.
Moreover, we assume that $\supp (f)\subset
\mathbb{N}^n$ and $f$ has no factor of the form $z^{\alpha}$ with
$\alpha =(\alpha_1,\ldots , \alpha_n)$.

\vspace{0.3cm}

The {\em amoeba}  $\mathscr{A}_f$ of an algebraic hypersurface $V_f\subset
 (\mathbb{C}^*)^n$
 is by definition ( see M. Gelfand, M.M. Kapranov
 and A.V. Zelevinsky \cite{GKZ-94}) the image of $V_f$ under the map :
\[
\begin{array}{ccccl}
\Log&:&(\mathbb{C}^*)^n&\longrightarrow&\mathbb{R}^n\\
&&(z_1,\ldots ,z_n)&\longmapsto&(\log\mid z_1\mid ,\ldots ,\log\mid
z_n\mid ).
\end{array}
\]

It was shown by M. Forsberg, M. Passare and A. Tsikh in \cite{FPT-00} that
there is an injective map between the set of components
$\{E^{c}\}$ of $\mathbb{R}^n\setminus \mathscr{A}_f$ and
$\mathbb{Z}^n\cap\Delta_f$:
$$
\ord_f :\{E^{c}\} \hookrightarrow \mathbb{Z}^n\cap\Delta_f
$$

\begin{The}[Forsberg-Passare-Tsikh, (2000)] Each  component
  of $\mathbb{R}^n\setminus \mathscr{A}_f$  is a convex domain and there
  exists a locally constant function:
$$
\ord_f :\mathbb{R}^n\setminus \mathscr{A}_f \longrightarrow \mathbb{Z}^n\cap\Delta_f
$$
which maps different components of the complement of $\mathscr{A}_f$
to different lattice points of $\Delta_f$.
\end{The}

\vspace{0.2cm}

Let $\mathbb{K}$ be the field of the Puiseux series with real power, which is the field of the  series $\displaystyle{a(t) = \sum_{j\in A_a}\xi_jt^j}$ with $\xi_j\in \mathbb{C}^*$ and $A_a\subset \mathbb{R}$ is a well-ordered set (which means that any subset has a smallest element). It is well known that the field $\mathbb{K}$ is algebraically closed and of characteristic zero, and it has a non-Archimedean valuation $\val (a) = - \min A_a$:
\[ 
\left\{ \begin{array}{ccc}
\val (ab)&=& \val (a) + \val (b) \\
\val (a + b)& \leq& \max \{ \val (a)  ,\, \val (b)  \} ,
\end{array}
\right.
\]
and we put $\val (0) = -\infty$. Let $g\in \mathbb{K}[z_1\ldots ,z_n]$ be a polynomial as in $(1)$ but the coefficients and the components of $z$ are in $\mathbb{K}$. If $<,>$ denotes the
scalar product in $\mathbb{R}^n$, then the following  piecewise affine linear convex function $\displaystyle{g_{trop} = \max_{\alpha\in\supp (g) } \{ \val (a_{\alpha}) + <\alpha , x>  \}}$, 
which is in the same time the Legendre transform of the function $\nu :\supp (g)\rightarrow \mathbb{R}$ defined by $\nu (\alpha ) = \min A_{a_{\alpha}}$, 
is called the {\em tropical polynomial} associated to $g$.

\begin{Def} The tropical hypersurface $\Gamma_g$ is the set of points in $\mathbb{R}^n$ where the tropical polynomial  $g_{trop}$ is not smooth (called the corner locus of  $g_{trop}$).
\end{Def}

\noindent We have the following Kapranov's theorem (see \cite{K-00}):

\begin{The}[Kapranov, (2000)] The tropical hypersurface $\Gamma_g$ defined by the tropical polynomial $g_{trop}$ is the subset of $\mathbb{R}^n$ image under  the valuation map of the algebraic hypersurface  defined by $g$. 
\end{The}

\noindent $\Gamma_g$ is also called the non-Archimedean amoeba of the zero locus of $g$ in $(\mathbb{K}^*)^n$.

\vspace{0.1cm} 

Let $g$ be a polynomial as above, $\Delta$ its Newton polytope, and $\tilde{\Delta}$ its  extending Newton polytope, i.e., $\tilde{\Delta} := convexhull \{ (\alpha , r )\in \supp (g)\times \mathbb{R} \mid \, r\geq 
\min A_{a_{\alpha}} \}$. Let us extend the above function $\nu$ (defined on $\supp (g)$) to all $\Delta$
 as follow:
\[
\begin{array}{ccccl}
\nu&:&\Delta&\longrightarrow&\mathbb{R}\\
&&\alpha&\longmapsto&\min \{ r\mid\, (\alpha ,r)\in \tilde{\Delta} \}.
\end{array}
\]
It's clear that the linearity domains of $\nu$ define a convex subdivision $\tau = \{\Delta_1,\ldots ,\Delta_l\}$ of $\Delta$ (by taking the linear subsets of the lower boundary of $\tilde{\Delta}$, see \cite{R-01}, \cite{PR1-04}, \cite{RST-05}, and \cite{IMS-07} for more details). Let $y= <x,v_i>+r_i$ be the equation of the hyperplane $Q_i\subset \mathbb{R}^n\times\mathbb{R}$ containing the points of coordinates $(\alpha ,\nu (\alpha ))$ with $\alpha \in \Verte (\Delta_i)$.

\noindent There is a duality between the subdivision $\tau$ and the subdivision of $\mathbb{R}^n$ induced by $\Gamma_g$ (see \cite{R-01}, \cite{PR1-04}, \cite{RST-05}, and \cite{IMS-07}), where each connected component of $\mathbb{R}^n\setminus \Gamma_g$ is dual to some vertex of $\tau_f$ and each $k$-cell of $\Gamma_g$ is dual to some $(n-k)$-cell of $\tau$. In particular, each $(n-1)$-cell of $\Gamma_g$ is dual to some edge of $\tau$.  If $x\in E_{\alpha\beta}^*\subset \Gamma_g$,  then $<\alpha , x> -\nu (\alpha ) = <\beta , x> -\nu (\beta )$, so $<\alpha  -\beta , x - v_i>  = 0$. This means that $v_i$ is a vertex of $\Gamma_g$ dual to some $\Delta_i$ having $E_{\alpha\beta}$ as edge.

\vspace{0.2cm}
\noindent More generally, let $\mathcal{I}\subset \mathbb{K}[z_1,\ldots ,z_n]$ be an ideal, and $V(\mathcal{I}) = \{ z\in(\mathbb{K}^*)^n \mid\, f(z)=0 \, for \,\,all \, f\in\mathcal{I}\}$. Let 
$x\in\mathbb{R}^n$, and $f\in \mathcal{I}$. The {\em initial form} $in_x(f)$ of $f$ is the polynomial in $\mathbb{C}[z_1,\ldots ,z_n]$ whose monomials are those which dominate $f$
in $\Val^{-1}(x)$. We denote by $in_x(\mathcal{I}) = \{ in_x(f)  \mid \, f\in \mathcal{I}\}$ the initial form of the ideal $\mathcal{I}$. It was shown by D. Speyer and B. Sturmfels that a tropical variety $\Val (V(\mathcal{I}))$ can be seen in several ways. A short and nice proof of this Theorem can be found in \cite{SS-04}:

\begin{The}[Speyer-Sturmfels, (2003)] For an ideal $\mathcal{I}\subset \mathbb{K}[z_1,\ldots ,z_n]$ the following subsets of $\mathbb{R}^n$ coincide:
\begin{itemize}
\item[(a)]\, The closure of the set $\{ (\val (z_1),\ldots ,\val (z_n)) \mid\, (z_1,\ldots ,z_n)\in V(\mathcal{I})\}$;
\item[(b)]\, The intersection of the tropical hypersurfaces defined by the tropical polynomials $f_{trop}$ where $f\in \mathcal{I}$;
\item[(c)]\, The set of all vectors $x\in \mathbb{R}^n$ such that $in_x(\mathcal{I})$ contains no monomial.
\end{itemize}
\end{The}

\vspace{0.2cm}

{\bf Passare-Rullg\aa rd function.}

\vspace{0.3cm}

Let $A'$ be the subset of $\mathbb{Z}^n\cap\Delta_f$ which is  the   image of
$\{E^{c}\}$ the set of complement components of the amoeba under the order mapping $\ord_f$. M. Passare and
H. Rullg\aa rd define in \cite{PR1-04}  the spine $\Gamma$ of the amoeba
$\mathscr{A}_f$  as the non-Archimedean amoeba defined by the
tropical polynomial:
$$
f_{trop}(x) = \max_{\alpha\in A'}\{ c_{\alpha}+<\alpha , x>\} ,
$$
For $\alpha \in A'$,   $c_{\alpha}$ is defined  by:
$$
c_{\alpha} = \R \left(   \frac{1}{(2\pi i)^n}\int_{\Log^{-1}(x)}\log
  \left(\frac{f(z)}{z^{\alpha}}\right) \frac{dz_1\wedge \ldots \wedge
    dz_n}{z_1\ldots z_n}\right)
$$
where  $x$ is in the complement component $E_{\alpha}^c$ of order $\alpha$,\, $z = (z_1,\cdots ,z_n)\in (\mathbb{C}^*)^n$. In other words, the spine of
$\mathscr{A}_f$ is defined as the set of points in $\mathbb{R}^n$
where the piecewise affine linear function $f_{trop}$ is not
differentiable, or as the projection in $\mathbb{R}^n$ of the corner locus of the 
graph of this function where $\mathbb{R}$ is
the semi-field $(\mathbb{R}; \max , +)$.
Let us denote by $\tau$ the convex subdivision of $\Delta_f$ dual
to the tropical variety $\Gamma$.

\vspace{0.2cm}

\begin{Def}{\rm (Ronkin ).} If $f$ is a Laurent polynomial, the 
function $\mathscr{N}_f$ defined in $\mathbb{R}^n$ by:
$$
\mathscr{N}_f(x)=\frac{1}{(2\pi i)^n}\int_{\log^{-1}(x)}\frac{\log\mid 
f(z)\mid dz_1\wedge\ldots\wedge dz_n}{z_1\ldots z_n}
$$
is called the Ronkin function of $f$.
\end{Def}
This function is convex and affine linear in any open connected set $E$
which is contained in $\mathbb{R}^n\setminus \mathscr{A}_f$. In
particular, if $f$ is a monomial, then $\mathscr{N}_f$ is an
affine linear function (i.e., if  $f(z,w)=a_\alpha z^\alpha$, then
$\mathscr{N}_f(x)= <\alpha ,x> + \log\mid a_\alpha\mid$). If $x\in \mathbb{R}^n\setminus \mathscr{A}_f$, the  
 gradient  of $\mathscr{N}_f$ on $x$  is the order 
 of the complement component containing $x$ 
 (\cite{FPT-00}).

\vspace{0.2cm}

\noindent
 We define the Passare-Rullg\aa rd's function $\nu_{PR} :\Delta_f\longrightarrow \mathbb{R}$
 on the Newton
polytope $\Delta_f$ as follows :
\begin{itemize}
\item[(i)]\, if $\alpha\in \Verte (\tau )$, then  $\nu_{PR} (\alpha
  )= -c_{\alpha}$
\item[(ii)]\, Let
  $\Delta_i$ be an element of the subdivision $\tau$ with maximal
  dimension, and  $y=<x,a_i>+b_i$ be the equation of the hyperplane in
  $\mathbb{R}^n\times \mathbb{R}$ containing the points of coordinates
  $(\alpha , -c_{\alpha})\in\mathbb{R}^n\times \mathbb{R}$ for
  $\alpha\in \Verte (\Delta_i )$, \, $a_i=(a_{1,\, i},\ldots ,a_{
  n,\, i})\in \mathbb{R}^n$ and $b_i\in\mathbb{R}$. 
If $\alpha\in\Delta_i\setminus \Verte (\tau )$, and for any small perturbation of the coefficient $a_{\alpha}$ the amoeba of the new polynomial remains without complement component of order $\alpha$, then   
we set $\nu_{PR} (\alpha ) = <\alpha ,a_i>+b_i +1$.Otherwise we set $\nu_{PR} (\alpha ) = <\alpha ,a_i>+b_i$.
\end{itemize}

If $f$ is the polynomial given by $(1)$, we define a family
of polynomials $\{ f_t\}_{t\in ]0,\frac{1}{e}]}$ as follows :
$$
f_t(z) = \sum_{\alpha\in \supp (f)} \xi_{\alpha}t^{\nu_{PR}
(\alpha )}z^{\alpha}
$$
where $\xi_{\alpha} = a_{\alpha} e^{\nu_{PR} (\alpha )}$ (note that $f_{\frac{1}{e}} = f$).

\begin{Rem} Passare and Rullg\aa rd prove that if $F$ is a face of the Newton polytope of $f$, and  $\alpha\in F$ is contained in the image of the order map, then $c_{\alpha}(f) = c_{\alpha}(f^{F})$ where $f^{F}$ is the truncation of $f$ to $F$. In particular, if $\alpha$ is a vertex of the Newton polytope of $f$, then $c_{\alpha} = \log \mid a_{\alpha}\mid$ (see \cite{PR1-04} Proposition 2).
\end{Rem}

\vspace{0.2cm}

\section{Complex tropical hypersurfaces}

\vspace{0.2cm}

 Let $h$ be  a strictly positive real number and $H_h$ be the self
 diffeomorphism of $(\mathbb{C}^*)^n$ defined by :
\[
\begin{array}{ccccl}
H_h&:&(\mathbb{C}^*)^n&\longrightarrow&(\mathbb{C}^*)^n\\
&&(z_1,\ldots ,z_n)&\longmapsto&(\mid z_1\mid^h\frac{z_1}{\mid
  z_1\mid},\ldots ,\mid z_n\mid^h\frac{z_n}{\mid z_n\mid} ). 
\end{array}
\]
which defines a new complex structure on $(\mathbb{C}^*)^n$
denoted by $J_h = (dH_h)\circ J\circ (dH_h)^{-1}$ where $J$ is the
standard complex structure.

\noindent A $J_h$-holomorphic hypersurface $V_h$ is a hypersurface
holomorphic with respect to the $J_h$ complex structure on
$(\mathbb{C}^*)^n$. It is equivalent to say that $V_h = H_h(V)$ where
$V\subset (\mathbb{C}^*)^n$ is an holomorphic hypersurface for the
standard complex structure $J$ on $(\mathbb{C}^*)^n$.

Recall that the Hausdorff distance between two closed subsets $A,
B$ of a metric space $(E, d)$ is defined by:
$$
d_{\mathcal{H}}(A,B) = \max \{  \sup_{a\in A}d(a,B),\sup_{b\in B}d(A,b)\}.
$$
Here we take $E =\mathbb{R}^n\times (S^1)^n$, with the distance
defined as the product of the
Euclidean metric on $\mathbb{R}^n$ and the flat metric on $(S^1)^n$.

\begin{Def} A complex tropical hypersurface $V_{\infty}\subset
  (\mathbb{C}^*)^n$ is the limit (with respect to the Hausdorff
  metric on compact sets in $(\mathbb{C}^*)^n$) of a  sequence of a
  $J_h$-holomorphic hypersurfaces $V_h\subset (\mathbb{C}^*)^n$ when
  $h$ tends to zero.
\end{Def}

\vspace{0.1cm}

\noindent  The argument map is the map  defined as follow:
\[
\begin{array}{ccccl}
\widetilde{\Arg}&:&(\mathbb{C}^*)^n&\longrightarrow&(S^1)^n\\
&&(z_1,\ldots ,z_n)&\longmapsto&(\widetilde{\arg} (z_1),\ldots ,\widetilde{\arg} (
z_n) ).
\end{array}
\]
We use the following notations: if $z = (z_1,z_2,\ldots ,z_n)\in (\mathbb{C}^*)^n$ and $z_j = \rho_je^{i\gamma_j}$, then $\widetilde{\Arg} (z) = (\widetilde{\arg} (z_1),\widetilde{\arg} (z_2),\ldots ,\widetilde{\arg} (z_n)) := (e^{i\gamma_1},e^{i\gamma_2},\ldots ,e^{i\gamma_n})$ and $\Arg (z) = (\arg (z_1),\arg (z_2),\ldots ,\arg (z_n)) := (\gamma_1,\gamma_2,\ldots ,\gamma_n)$.

\vspace{0.1cm}

\noindent If we denote by $\tau_{\tilde{a}}$ the translation in $\mathbb{R}^n$ by the vector ${}^tL^{-1}(\arg (a_1),\ldots ,\arg (a_n))$, then we have the following commutative diagram:
\begin{equation}
\xymatrix{
V_f\ar[d]_{\Arg}\ar[rr]^{\Phi_{L,\, a}}&&P_1\ar[d]^{\Arg}\cr
\mathbb{R}^n\ar[rr]^{{}^tL\circ\tau_{\tilde{a}}}&&\mathbb{R}^n.
}\nonumber
\end{equation}
We have the same diagram if we replace $\Arg$ by the logarithmic map $\Log$
  (or the valuation map $\Val$ if we work in $(\mathbb{K}^*)^n$).
\vspace{0.1cm}
\noindent We complexify the valuation
map as follows :
\[
\begin{array}{ccccl}
w&:&\mathbb{K}^*&\longrightarrow&\mathbb{C}^*\\
&&a&\longmapsto&w(a ) = e^{\val (a )+i{\arg} (\xi_{-\val
    (a )})}
\end{array}
\]
Let  $\widetilde{\Arg}$ be the argument map $\mathbb{K}^*\rightarrow S^1$ defined by: for any $ a \in \mathbb{K}$ with 
$\displaystyle{a = \sum_{j\in A_a}\xi_jt^j}$, $\widetilde{\Arg} (a) = e^{i\arg (\xi_{-\val (a)})}$ (this map extends the map $\widetilde{\Arg} : \mathbb{C}^*\rightarrow S^1$ defined by $\rho e^{i\theta} \mapsto e^{i\theta}$).

Applying this map coordinate-wise we obtain a map :
\[
\begin{array}{ccccl}
W:&(\mathbb{K}^*)^n&\longrightarrow&(\mathbb{C}^*)^n
\end{array}
\]

Using Kapranov's theorem \cite{K-00} and degeneration of a complex structures, 
Mikhalkin gives an algebraic
definition of a complex tropical hypersurfaces (see \cite{M2-04}) as
follows:

\begin{The}[Mikhalkin, (2002)]  The set  $V_{\infty}\subset (\mathbb{C}^*)^n$
  is a complex tropical hypersurface if and only if there
  exists an algebraic hypersurface
  $V_{\mathbb{K}}\subset(\mathbb{K}^*)^n$ such that $W(V_{\mathbb{K}}) = V_{\infty}$.
\end{The}

Let $\Log_{\mathbb{K}}(z_1,\ldots ,z_n)=(\val (z_1),\ldots ,
\val (z_n))$, which means that $\mathbb{K}$ is equipped with the
norm defined by $\Vert z\Vert_{\mathbb{K}}= e^{\val (z)}$ for any
$z\in\mathbb{K}^*$.
Then we have the following commutative diagram:
\begin{equation}
\xymatrix{
(\mathbb{K}^*)^n\ar[rr]^{W}\ar[dr]_{\Log_{\mathbb{K}}}&&(\mathbb{C}^*)^n\ar[dl]^{\Log}\cr
&\mathbb{R}^n
}\nonumber
\end{equation}

\vspace{0.2cm}

\subsection[]{Complex tropical hypersurfaces with a simplex Newton polytope} ${}$

\vspace{0.2cm}

Let $K$ be the field of the complex number $\mathbb{C}$ or the field of the generalized Puiseux series $\mathbb{K}$.
Let $a=(a_1,\ldots ,a_n)\in (K^*)^n $ and  $P_a\subset(K^*)^n$ be the hyperplane defined by the polynomial $f_a(z_1,\ldots ,z_n) = 1+\sum_{j=1}^na_jz_j$, and  $P_1$ the one defined by $f_1(z_1,\ldots ,z_n) =1+\sum_{j=1}^nz_j $, then it's clear that 
if $\tau_{a^{-1}}$ is the translation in the multiplicative group $(K^*)^n$ by $a^{-1}$, then  we have
$P_a = \tau_{a^{-1}}(P_{1})$. Let $L$ be an invertible matrix with integer coefficients and positive determinant
$$
L = \left(
\begin{array}{ccc}
\alpha_{11}&\hdots&\alpha_{1n}\\
\vdots&\ddots&\vdots\\
\alpha_{n1}&\hdots&\alpha_{nn}
\end{array}
\right) ,
$$
and $\Phi_{L,\, a}$ be the homomorphism of the algebraic torus defined as follow.
\[
\begin{array}{ccccl}
\Phi_{L,\, a}&:&(K^*)^n&\longrightarrow&(K^*)^n\\
&&(z_1,\ldots ,z_n)&\longmapsto&(a_1\prod_{j=1}^nz_j^{\alpha_{1j}},\ldots ,a_n\prod_{j=1}^nz_j^{\alpha_{nj}}).
\end{array}
\]
Let $V_f\subset (\mathbb{K}^*)^n$ be the hypersurface defined by the polynomial 
$$
f(z_1,\ldots ,z_n)=1+\sum_{k=1}^na_k\prod_{j=1}^nz_j^{\alpha_{kj}}.
$$
Let us denote by 
$co\mathscr{A}(V_f)$ 
the set of arguments of $V_f$, and by the same $co\mathscr{A}(V_f)$ its lifting in the universal covering $\mathbb{R}^n$ of the real torus (abuse of notation). Let $\tau_{\Arg (a)}$ be the translation by the vector $(\arg (a_1),\ldots ,\arg (a_n))$ in $\mathbb{R}^n$.

\begin{Lem} Let $V_f$ be a hypersurface defined by a polynomial $f$ with Newton polytope a simplex $\Delta$, and such that its support $\supp (f)$ is precisely the vertices of $\Delta$. Then we have:
\begin{itemize}
\item[(i)]\, if $K = \mathbb{C}$ (resp. $\mathbb{K}$), then the amoeba $\mathscr{A}_f$ of the hypersurface $V_f$ is the image under ${}^tL^{-1}\circ \tau^{-1}_{\Log (a)}$ (resp. ${}^tL^{-1}\circ \tau^{-1}_{\Val (a)}$) of the amoeba of the standard hyperplane $P_1$. In particular it is solid.
\item[(ii)]\, if $K = \mathbb{C}$  or $\mathbb{K}$, then the coamoeba $co\mathscr{A}_f$ of the hypersurface $V_f$ is the image under ${}^tL^{-1}\circ \tau^{-1}_{\Arg (a)}$ of the coamoeba of the standard hyperplane $P_1$. In particular, the number of its complement components in the real torus $(S^1)^n$ is equal to $n!\Vol (\Delta)$.
\end{itemize}
\end{Lem}

This means that we have the following:
$$
\hspace{0.3cm} {}^tL^{-1}\circ \tau^{-1}_{\Log (a)} (\mathscr{A}_{P_1}) = \mathscr{A}_f, \quad\mbox{and}\quad
{}^tL^{-1}\circ \tau^{-1}_{\Arg (a)} (co\mathscr{A}_{P_1}) = co\mathscr{A}_f. 
$$

\begin{prooof} First of all, we can see that the Newton polytope $\Delta_f$ of $f$ is the image under the linear map $L$ of the standard simplex. The matrix $L$ is invertible, so $\Phi_{L,\, a}(V_f) = P_1$. Indeed, if $(z_1',\ldots ,z_n')$ is in $P_1$, then there exists $(u_1,\ldots ,u_n)\in (\mathbb{C}^*)^n$ such that for any $1\leq j\leq n$, we have $z_j'= a_je^{u_j}$. The matrix $L$ is invertible, so its column vectors $\alpha_k$ are linearly independent. Hence, there exists a vector $(v_1,\ldots ,v_n)\in (\mathbb{C}^*)^n$ which is a solution of the following linear system:
$$
\begin{array}{ccccccc}
\alpha_{11}x_1&+&\hdots&+&\alpha_{1n}x_n&=& u_1\\
\alpha_{21}x_1&+&\hdots&+&\alpha_{2n}x_n&=& u_2\\
\vdots&\vdots&\vdots&\vdots&\vdots&\vdots&\vdots\\
\alpha_{n1}x_1&+&\hdots&+&\alpha_{nn}x_n&=& u_n;
\end{array}
$$
and then 
$$
\begin{array}{lll}
\Phi_{L,\, a} (e^{v_1},\ldots ,e^{v_n}) &= & (a_1\prod_j e^{\alpha_{1j}v_j},\ldots ,a_n\prod_j e^{\alpha_{nj}v_j})\\
&=& (a_1e^{u_1},\ldots ,a_ne^{u_n})\\
&=& (z_1',\ldots ,z_n').
\end{array}
$$
But  $(e^{v_1},\ldots ,e^{v_n})\in V_f$, because $1+ \sum_k a_k\prod_j (e^{v_j})^{\alpha_{kj}} = 1 + \sum_k a_k e^{\sum_jv_j\alpha_{kj}}  = 1 + \sum_k a_k e^{u_k} = 1 + \sum_k z_k' = 0$. So, $\Phi_{L,\, a}(V_f) = P_1$, and the Lemma is done after using the properties of the logarithmic and the argument maps on one hand, and the properties of the amoeba and the coamoeba of the standard hyperplane on the other hand, and the fact that $\det (L) = n!\Vol (\Delta )$. Recall that the amoeba of the standard hyperplane is solid because of the injectivity of the order map, and the fact that the standard simplex contains no interior lattice point. 
\end{prooof}

\noindent We have the same properties   for the complex tropical hypersurface $W(V_{\mathbb{K}})$ i.e., ${}^tL^{-1}\circ \tau^{-1}_{\Arg (a)} (co\mathscr{A}_{P_1}) = co\mathscr{A}_{W(V_{\mathbb{K}})}$. 
Therefore, the set of arguments of  any hypersurface defined by a maximally sparse polynomial with Newton polytope a simplex can be easily drawn.

\vspace{0.1cm}

\begin{Exe}

We draw in figure 1 the coamoeba  of the complex  curve defined by the polynomial $f_1(z_1,z_2)= 1+ z_1^2z_2^3+z_1^3z_2$
where the matrix  ${}^tL_1^{-1}$ is equal to $\frac{1}{7}\left(
  \begin{array}{cc} 3&-1\\ -2&3\end{array}\right) $
and in figure 2 the coamoeba of the complex curve 
defined by the polynomial $f_2(z_1,z_2) = z_1^2z_2^2+z_1+z_2$
where the matrix  ${}^tL_2^{-1}$ is equal to $\frac{1}{3}\left(
  \begin{array}{cc} 1&1\\ -2&1\end{array}\right)$.

\begin{figure}[h!]
\begin{center}
\includegraphics[angle=0,width=0.4\textwidth]{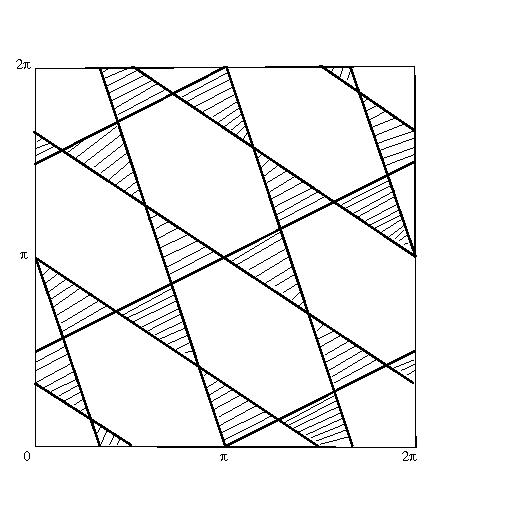}
\caption{The coamoeba of the curve defined by the polynomial $f_1$}
\label{c}
\end{center}
\end{figure}

\begin{figure}[h!]
\begin{center}
\includegraphics[angle=0,width=0.4\textwidth]{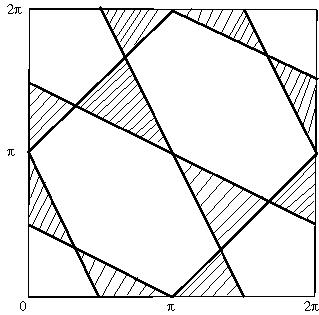}
\caption{The coamoeba of the curve defined by the polynomial $f_2$}
\label{c}
\end{center}
\end{figure}

\end{Exe}

\section{Codual hyperplanes arrangement and Non-Archimedean Coamoebas}

\vspace{0.2cm}

Let $A$ be a finite subset in $\mathbb{Z}^n$, and $\nu_1$ be a convex function on $A$. Let us denote by $\Delta$ the convex hull of $A$. Moreover, assume that the extension of $\nu_1$ to $\Delta$ defines a convex triangulation $\tau$  of $\Delta$ such that $\Verte (\tau ) = A$. Let $\nu$ be the map defined on $A$ which take its values in $\mathbb{C}^*$, such function is called the {\em generalized Passare-Rullg\aa rd function}:  
\[
\begin{array}{ccccl}
\nu&:&A&\longrightarrow&\mathbb{C}^*\cong \mathbb{R}\times S^1\\
&&\alpha&\longmapsto&(\nu_1(\alpha ), \nu_2(\alpha )).
\end{array}
\]
such that for any real positive number $t$, the number  $t^{-\nu_1(\alpha )}\nu_2(\alpha )$ can be seen as the complex number of argument
$\nu_2(\alpha )$ and of norm $t^{-\nu_1(\alpha )}$. 
Let $f$ be the polynomial in $\mathbb{K}[z_1,\ldots ,z_n]$ defined by:
$$
f(z) = \sum_{\alpha\in A}\nu_2(\alpha )t^{\nu_1(\alpha )}z^{\alpha},
$$
and $V$ be the algebraic hypersurface over $\mathbb{K}$ with defining polynomial $f$. Let $V_{\infty} = W(V)$ be the complex tropical hypersurface associated to $f$.
For any edge $E_{\alpha\beta}$ in $\tau$ with extremities $\alpha$ and $\beta$, let us denote by $H_{\alpha\beta}$ its {\em codual hyperplane} relatively to $\nu$ i.e., the hyperplane of arguments of the holomorphic cylinder defined by the polynomial $\nu (\alpha )z^{\alpha} + \nu (\beta )z^{\beta} = 0$:
$$
H_{\alpha\beta} := \Arg \{ z\in (\mathbb{C}^*)^n \mid\, \nu (\alpha )z^{\alpha} + \nu (\beta )z^{\beta} = 0\};
$$
if we denote by $x = (\arg (z_1),\ldots ,\arg (z_n))$, then its equation is given by:
$$
\arg (\nu_2 (\alpha )) - \arg (\nu_2 (\beta )) + <\alpha -\beta , x> = \pi + 2k\pi
$$
with $k\in \mathbb{Z}$.

\vspace{0.1cm}

\begin{Def} The coamoeba of complex tropical variety is called {\em Non-Archimedean coamoeba}.
\end{Def}

\vspace{0.1cm}

\begin{Rems} ${}$

\begin{itemize}
\item[(1)]\, The hyperplanes $H_{\alpha\beta}$ can be seen as a real  $(n-1)$-torus in $(S^1)^n$, or as the image of a hyperplane of rational slope in $\mathbb{R}^n$ by the group of translations $(2\pi\mathbb{Z})^n$.
\item[(2)]\, If $E_{\alpha\beta}$ is an external edge in $\tau$, i.e., an edge of the $\Delta$, then $H_{\alpha\beta}$ is called an {\em external codual hyperplane}. In this case $H_{\alpha\beta}$ is equipped with a frame for any element $\Delta_v$  of $\tau$ containing $E_{\alpha\beta}$ as edge by the following: the edge $E_{\alpha\beta}$ is adjacent to $n-1$ facets of $\Delta_v$, we take a primitive outward normal vector for each one, this gives a frame of $H_{\alpha\beta}$ after taking an order of those normal vectors such that they give a direct base of $H_{\alpha\beta}$.
\item[(3)]\, We can associated to any codual hyperplane $H_{\alpha\beta}$ a natural number different than zero as follow: $w(H_{\alpha\beta})$ is the cardinality of edges in $\tau$ such that their codual hyperplane is $H_{\alpha\beta}$, this number is called the {\em weight} of $H_{\alpha\beta}$.
\end{itemize}
\end{Rems}

\begin{Def} Let $\nu$ be a map as above. A {\em weighted  arrangement} $\mathscr{H}_{\nu}$  of codual hyperplanes relatively to $\nu$  of a triangulation $\tau$ defined by $\nu_1$ is the union of all the codual hyperplanes $H_{\alpha\beta}$ to the edges $E_{\alpha\beta}$ of $\tau$ relatively to $\nu$ equipped with their weight. With the notation as above, we say that $\mathscr{H}_{\nu}$ is the {\em codual hyperplanes arrangement associated to} $V$ (or $V_{\infty}$).
\end{Def}

Let  $\tau  = \{ \Delta_1,\ldots ,\Delta_l\}$ and for $i=1,\dots ,l$, let $L_i$ be the linear part of  a surjection (affine linear) of $\mathbb{R}^n$ which sends the standard simplex $\Delta_{std}$ to $\Delta_i$, which we assume of positive determinant. 

\begin{Rems} ${}$

\begin{itemize}
\item[(1)]\, Let $\mathscr{H}_{std}$ be the arrangement of the codual hyperplanes to the edges of the standard simplex $\Delta_{std}$ (I mean the Newton polytope of the hyperplane in $(\mathbb{K}^*)^n$ or in  $(\mathbb{C}^*)^n$ defined by the polynomial $f(z) = 1 + \sum_{j=1}^n z_j$). For  any $\Delta_i\in \tau$,  the arrangement $\mathscr{H}_{\Delta_i}$ of codual hyperplanes to the edges of $\Delta_i$ is equal to ${}^tL_i^{-1}(\mathscr{H}_{std})$.
\item[(2)]\, Let $E$ be a common edge to $\Delta_i$ and $\Delta_j$. We denote by  $E_{std}$ and $E_{std}'$ the edges of the standard simplex such  that $E = L_i(E_{std}) = L_j(E_{std}')$. Then we have ${}^tL_i^{-1}({H}_{std}) = {}^tL_j^{-1}({H}_{std}')$, where ${H}_{std}$ (resp. ${H}_{std}'$) is the codual hyperplane to ${E}_{std}$ (resp. ${E}_{std}'$).
\end{itemize}
\end{Rems}

\begin{figure}[h!]
\begin{center}
\includegraphics[angle=0,width=0.3\textwidth]{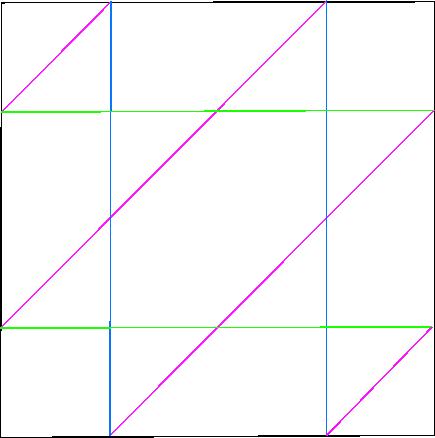}
\caption{Arrangement of the codual lines to the edges of the Newton polygon of the standard line; we draw here only four fundamental domain in the torus; the  standard line in $(\mathbb{C}^*)^2$ is defined by the polynomial $f(z,w)=z+w+1$.}
\label{}
\end{center}
\end{figure}

\begin{figure}[h!]
\begin{center}
\includegraphics[angle=0,width=0.3\textwidth]{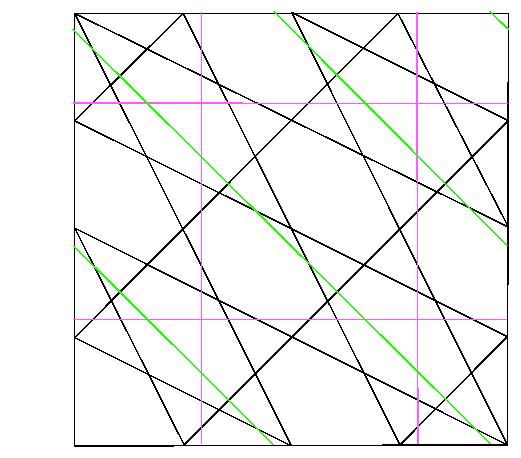}\quad
\includegraphics[angle=0,width=0.3\textwidth]{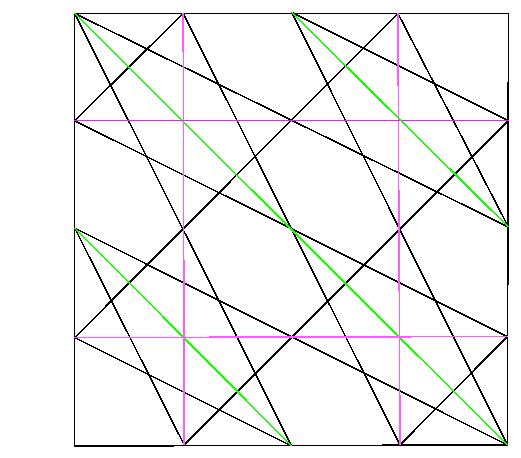}\quad
\includegraphics[angle=0,width=0.2\textwidth]{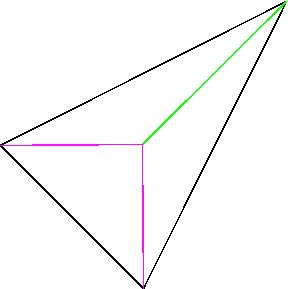}
\caption{Two arrangements of the codual lines to the edges of the triangulation  of the Newton polygon  in the right with two different combinatorial types; we draw here only four fundamental domains in the torus $(S^1)^2$.}
\label{}
\end{center}
\end{figure}

\begin{figure}[h!]
\begin{center}
\includegraphics[angle=0,width=0.3\textwidth]{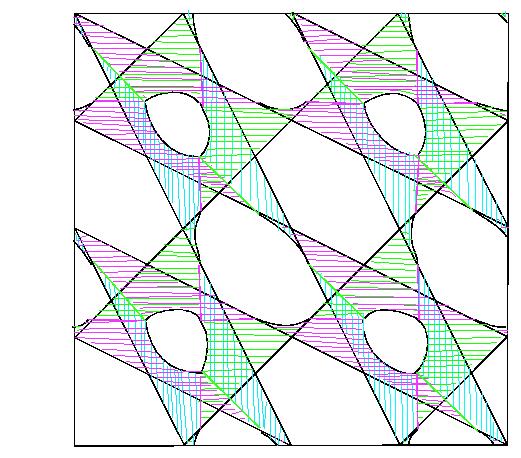}\quad
\includegraphics[angle=0,width=0.3\textwidth]{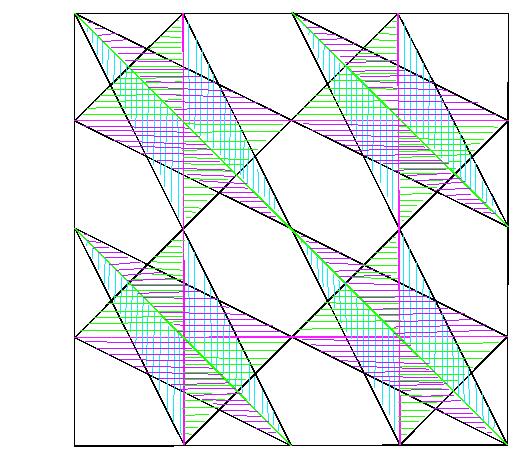}
\caption{Their associated complex coamoebas}
\label{}
\end{center}
\end{figure}

\begin{Def} A perturbation of an arrangement $\mathscr{H}$ of codual hyperplanes is called {\em small} if there is a point  $p$ equal to  the intersection of $k$ codual hyperplanes in $\mathscr{H}$  with $k\geq n$, then the image of $p$ under the perturbation is also the intersection of the image of the  $k$ codual hyperplanes under this perturbation.
\end{Def}

\begin{Def} Let $\tau$ be the convex subdivision of the  convex hull $\Delta$  of a finite set of  lattice points $A$ in $\mathbb{Z}^n$. Assume that $\tau$ is  defined by a strictly convex function $\nu_1 : A\rightarrow \mathbb{R}$. If  $\nu_2$ and $\nu_2'$  are two maps $A\rightarrow S^1 $,  and $\mathscr{H}_{\nu_2}$,\, $\mathscr{H}_{\nu_2'}$  are their associated  arrangement of codual hyperplanes of the edges of $\tau$,  we say that $\mathscr{H}_{\nu_2}$ and $\mathscr{H}_{\nu_2'}$ are in the same {\em combinatorial type} if  one can obtain $\mathscr{H}_{\nu_2'}$  by a  continuous small perturbation of $\mathscr{H}_{\nu_2'}$ in the set of codual hyperplanes
\end{Def}

\begin{Pro} With notations as above, the number of combinatorial type of arrangement of codual hyperplanes to the edges of $\tau$ is finite.
\end{Pro}

\begin{prooof} it suffices to prove the finiteness of combinatorial types for the arrangements with fixed values of  $\nu_2$ on the vertices of $\Delta$.
\end{prooof}

\begin{Def} Let $V$ be a complex algebraic hypersurface defined by a polynomial $f$ with Newton polytope $\Delta$. Let $\Gamma$ the spine of the  amoeba of $V$, and $\tau$ be its dual convex subdivision of $\Delta$. An edge $E_{\alpha\beta}$ of ends $\alpha$ and $\beta$ is called {\em effective} if $\alpha$ and $\beta$ are in the intersection of the  support of $f$ with the set of vertices of $\tau$; otherwise it is called {\em virtual}.
\end{Def}

\begin{Def} An arrangement $\mathscr{H}$ of codual hyperplanes to the edges of $\tau$ is called {\em effective} (resp.{\em virtual}) if any hyperplane $H$ in $\mathscr{H}$ is codual to an effective  (resp. virtual) edge of $\tau$.
\end{Def}

\begin{Def} Let  $a_{\alpha}$  be the coefficient of index $\alpha$ of a  complex polynomial $f$. Then the  coefficient $a_{\alpha}$  is said {\em virtually zero} (resp. {\em virtually  non zero}) if $a_{\alpha} = 0$ and  $\alpha\in\Ima (\ord )$ (resp. $a_{\alpha} \ne 0$ and  $\alpha\notin\Ima (\ord )$).
\end{Def}

\begin{Def} Let $a_{\alpha}$  be the coefficient of index $\alpha$ of a complex polynomial $f$. Then the coefficient $a_{\alpha}$  is said {\em effectively  zero} (resp. {\em effectively  non zero}) if $a_{\alpha} = 0$ and  $\alpha\notin\Ima (\ord )$ (resp. $a_{\alpha} \ne 0$ and  $\alpha\in\Ima (\ord )$).
\end{Def}

If the polynomial is non-Archimedean $\Ima (\ord )$ is replaced by  the set of $\alpha\in\Verte (\tau )$ where $\tau$ is  the dual subdivision of the tropical hypersurface  defined by the tropical polynomial associated to $f$.

\vspace{0.2cm}

\section{Geometry of coamoebas and their combinatorial structure}

\vspace{0.2cm}

Let $V$ be a complex algebraic hypersurface defined by the polynomial $f(z) = \sum_{\alpha\in\supp (f)}a_{\alpha}z^{\alpha}$.
 Let $\tau$ be the subdivision of $\Delta$ dual to the spine of the amoeba of $V$, and $c_{\alpha} = \R (\Phi_{\alpha}(f))$ is the constant of Passare-Rullg\aa rd defined in section 1 for any $\alpha\in \Ima (\ord_f)$. If $\Delta_i\in \tau$, then there exist $(a_{v_i}, b_{v_i})\in \mathbb{R}^n\times\mathbb{R}$ such that the hyperplane in $\mathbb{R}^n\times\mathbb{R}$ containing the extended $\Delta_i$ is given by the equation $y = <x , a_{v_i}> + b_{v_i} $.

\begin{Def} Let $\alpha$ be a lattice point in $\supp (f)\setminus \Ima (\ord_f)$. If there exists $c\in S^1$ such that for any $\varepsilon >0$, the perturbation of the polynomial $f$ given by replacing the coefficient $a_{\alpha}$ by $\varepsilon ca_{\alpha}$ gives a polynomial $f^{(\varepsilon)}$ such that $\alpha\in \Ima (\ord_{f^{(\varepsilon)}})$, then we say that $\alpha$ is of {\em type $I$}. Otherwise we say that $\alpha$ is of {\em type $II$}.
\end{Def}

\vspace{0.1cm}

\noindent The purpose of this section is to show that there exists a hyperplanes arrangement in the universal covering of the real torus which completely determines the  geometry and the topology of the coamoeba of $V$. Moreover, we prove that this arrangement has the same combinatorial type of a codual hyperplanes arrangement associated to a complex tropical hypersurface whose  image under the logarithmic map is the amoeba spine of $V$.

\vspace{0.1cm}

\begin{The} Let $V$ be complex hypersurface defined by a polynomial $f$ with Newton polytope $\Delta$. Let $A:=\Ima (\ord_f)\cup\{ type \,\, I\,\, indices\}$ be the union of the image of the order map and the indices in $\Delta$ of type I. Then there exists a function $\nu_{(2,\, f)} : A\rightarrow S^1$ such that
if $V_{\mathbb{K}}$ denotes the hypersurface over $\mathbb{K}$ defined by the polynomial:
$$
f_{\mathbb{K}}(z) = \sum_{\alpha\in A_f}b_{\alpha}z^{\alpha},
$$
with $b_{\alpha}(t)  = t^{\nu_{PR}(\alpha )}\nu_{(2,\, f)}(\alpha )$, then  the coamoeba of the complex tropical hypersurface $V_{\infty ,\, f} := W(V_{\mathbb{K}})$ is homeomorphic to the closure in the real torus of the coamoeba of $V$.
\end{The}

\vspace{0.1cm}

We can remark that the image under the logarithmic map of $W(V_{\mathbb{K}})$  is the amoeba spine of $V$.

\vspace{0.1cm}

\subsection{Definition and Geometric interpretation  of  the function $\nu_{(2,\, f)}$} ${}$

\vspace{0.1cm}

Let $V$ be a complex algebraic hypersurface with defining polynomial $f$ as before, and assume, for simplicity,  that $\Delta \cap \mathbb{Z}^n$ contains no point of type I. The goal of this section is the construction of a function on $A_f := \supp (f) \cup \Ima (\ord_f)$ taking its values in $S^1$, and an algebraic hypersurface $V_{\mathbb{K}}$ over the field of Puiseux series $\mathbb{K}$, such that the coamoeba of its associated complex tropical hypersurface $W(V_{\mathbb{K}})$ contains a codual hyperplanes arrangement $\mathscr{H}_{\nu}$ where $\nu = (\nu_{PR} ,\, \nu_{(2,\, f)} ) : A_f\rightarrow \mathbb{R}\times S^1$. Moreover, the defining polynomial of  $V_{\mathbb{K}}$ is given by:
$$
f_{\mathbb{K}}(z) = \sum_{\alpha\in A_f}b_{\alpha}z^{\alpha},
$$
with $b_{\alpha}(t)  = t^{\nu_{PR}(\alpha )}\nu_{(2,\, f)}(\alpha )$
and the image under the logarithmic map of $W(V_{\mathbb{K}})$  is the amoeba spine of $V$. We denote by $f_t$ the following complex polynomial:
$$
f_t(z) = \sum_{\alpha\in A_f}t^{\nu_{PR}(\alpha )}\nu_{(2,\, f)}(\alpha )z^{\alpha}.
$$ 

\vspace{0.2cm}

\noindent Before beginning the construction of $\nu_{(2,\, f)}$ for any $n$, let us  looking at the following example in one  variable:

\vspace{0.2cm}

\begin{Exe} Let $f$ be the complex polynomial in one variable defined by: 
$$
f(z) = (z-r_1)(z-r_2)\cdots (z-r_{N}) = z^N + a_{N-1}z^{N-1} + \cdots + a_0,
$$
and we denote by $V_f$ the zero locus of $f$. Assume that the roots of $f$ are such that $0 <\mid r_1\mid < \mid r_2\mid < \cdots <\mid r_N\mid$. Let $0 < k < N$, then up to a constant we have:
$$
\Arg (\nu_{(2,\, f)}(k )) =  (1 + (-1)^{N-k +1})\frac{\pi}{2} + \sum_{j=k +1}^N\Arg (r_j) \qquad\qquad \qquad (2).
$$
Indeed, in this case we have $f_{\mathbb{K}}(z) = \sum_{k = 0}^N t^{\nu_{PR}(k)}
\nu_{(2,\, f)}(k )z^{k}$
because all the lattice points in the interval $[0; N]$ are in the image of the order map. In particular, we have $H_{k ,\, k -1} := \Arg \left( \{ z\in \mathbb{C}^*\, \mid\, \nu_{(2,\, f)}(k )t^{\nu_{PR}(k )}z^{k} + \nu_{(2,\, f)}(k -1 )t^{\nu_{PR}(k -1 )}z^{k -1} = 0  \}\right)$. So, if we denote by $x$ the argument of $z$ we obtain $\Arg (\nu_{(2,\, f)}(k )) - \Arg (\nu_{(2,\, f)}(k -1)) + x = \pi \mod (2\pi) $. Let $k = N$, then we have:
$$
\Arg (\nu_{(2,\, f)}(N)) - \Arg (\nu_{(2,\, f)}(N -1)) + x = \pi \mod (2\pi).
$$
But $x = \Arg (r_N) \mod (2\pi)$, because if we look to the amoeba, this means that $\Log (\Arg^{-1}(x) \cap W(V_{\mathbb{K}})) = \log (r_N) = \log\mid r_N\mid$. Hence $\Arg (\nu_{(2,\, f)}(N -1)) =  \Arg (\nu_{(2,\, f)}(N)) + \Arg (r_N) + \pi \mod (2\pi)$.
If  $k = N-1$, we obtain :
$$
\Arg (\nu_{(2,\, f)}(N -1)) - \Arg (\nu_{(2,\, f)}(N -2)) + x = \pi \mod (2\pi),
$$
and  $x = \Arg (r_{N-1}) \mod (2\pi)$ by  the same reasoning as before. So
$$
\Arg (\nu_{(2,\, f)}(N -2)) = \Arg (\nu_{(2,\, f)}(N)) + \Arg (r_{N-1}) + \Arg (r_N) + \pi + \pi \mod (2\pi),
$$
Hence, for $0 < k < N$ we have the formula $(2)$, up to multiplication by $\nu_{(2,\, f)}(N)$.
\end{Exe}

\vspace{0.2cm}

\noindent {\bf  First step: construction and definition of $\nu_{(2,\, f)}$ for $n=1$.} ${}$

\vspace{0.2cm}

\noindent First of all, let us define $\nu_{(2,\, f)}$ on the set of vertices of $\Delta$, and on all the lattice points in the intersection of the  external edges of  $\Delta$ with $\Ima (\ord_f)$. Also, we give 
 the geometric  reasons for which this definition is justified. If  $\alpha\in\Verte (\Delta )$ then we put $\nu_{(2,\, f)}(\alpha )$. The justification is the following:
let  $E_{\alpha\beta}$ be an external edge of $\Delta$,  and consider the truncation of $f$ to $E_{\alpha\beta}$ as a polynomial $l$ in one variable. So, by multiplying $l$ by a complex constant we can assume that 
$$
l(z) = \prod_{j=1}^n\prod_{s=1}^{l_j}(z-r_je^{i\theta_{js}}) = z^N+a_{N-1}z^{N-1}+ \cdots +a_0
$$ 
with  $N = \sum_{j=1}^n l_j$, and the $r_j$'s are positive real numbers such that $r_1<r_2<\cdots <r_n$. The index $0 \leq k \leq n$ are the lattice  points in  $E_{\alpha\beta} \cap \Ima (\ord_f)$, where $\alpha$ is replaced by $0$ and $\beta$ by $n$.  We can remark that the amoeba of $l$ has $n$ points which are $\log (r_j)$. Let us define $\nu_{(2,\, f)}$ on $E_{\alpha\beta} \cap \Ima (\ord_f) $, and take $0<k<n$, such that  the following equation holds:
$$
\Arg (\nu_{(2,\, f)}(k)) - \Arg (\nu_{(2,\, f)}(k-1)) + l_k x = \pi \mod (2\pi ),\qquad\qquad (3)
$$
and if we replace $x$ in $(3)$ by the average of the arguments of all roots of module $r_k$, we obtain:
$$
\Arg (\nu_{(2,\, f)}(k)) - \Arg (\nu_{(2,\, f)}(k-1)) + l_k (\frac{\sum_{s=1}^{l_k} \theta_{ks}}{l_k}) = \pi \mod (2\pi ).
$$
Recall that we can put $\nu_{(2,\, f)} (n) = 1$, so we have:
$$
\Arg (\nu_{(2,\, f)} (k)) = (1+(-1)^{n-k+1})\frac{\pi}{2} + \sum_{j=k+1}^{n}\sum_{s=1}^{l_{k+1}}\theta_{(k+1)s}.
$$
We can remark that in this case $\Arg (\nu_{(2,\, f)} (\alpha ))$ is equal to $\Arg (a_{\alpha})$  for any $\alpha\in \Verte (\Delta )$. Indeed, if we don't normalize the polynomial $l$ then we obtain $\Arg (\nu_{(2,\, f)} (\alpha )) = \Arg (\nu_{(2,\, f)} (\beta )) + \Arg (a_{\alpha}) - \Arg (a_{\beta})$, which means that $\Arg (\nu_{(2,\, f)} (\alpha )) - \Arg (a_{\alpha})$ is constant and don't depends on $\alpha$ if this index is in the set of vertices of $\Delta$. Hence, by multiplying the original polynomial $f$ by a constant if necessary, we obtain $\Arg (\nu_{(2,\, f)} (\alpha )) = \Arg (a_{\alpha})$ for any $\alpha\in \Verte (\Delta )$.

\vspace{0.2cm}

\begin{Rems}${}$

\begin{itemize}
\item[(1)]\, It may be noted that for $n = 1$ , the codual hyperplanes arrangement defined in section 4 coincide with the coamoeba of the complex tropical hypersurface (all are point in this case). Moreover, it is an average of the coamoeba of  the complex variety defined by the polynomial $f$. Another thing, if $\Delta$ contains an index $k$ of type I, this means that the polynomial $f$  has a multiple root.
\item[(2)]\, In dimension $n > 1$, we have a similar geometric interpretation, i.e., the coamoeba of the complex tropical hypersurface defined by the polynomial $f_{\mathbb{K}}$ is an average of  the coamoeba of the complex hypersurface defined by the polynomial $f$. Moreover, an index $\beta\in\Delta$ is of type I, means that the codual hyperplanes to the adjacent edges to $\beta$ intersect at a point. We will prove these affirmations  this section. 
\item[(3)]\, Recall, that if an index $\beta\in\Delta$ is of type I, this means that the normal vectors to the hyperplanes in $\mathbb{R}^n\times\mathbb{R}$ containing the extended elements of $\tau$ adjacent to $\beta$ coincide. The coordinates of these normal vectors are those  of the vertices in the cycle of the non-Archimedean amoeba  bounding the complement component of order $\beta$. It is the geometric interpretation in terms of non-Archimedean amoebae.
\end{itemize}
\end{Rems}

\vspace{0.2cm}

\noindent {\bf  Second  step: definition  of $\nu_{(2,\, f)}$ for any $n$.} ${}$

\vspace{0.2cm}

\noindent Now, we are able to define $\nu_{(2,\, f)}$ on the lattice points in the intersection of $\Ima (\ord_f)$  with  any  $s$-cells of $\tau$ in the boundary of $\Delta$ for $0\leq s\leq n$. Let  $\alpha\in\Ima (\ord_f)$, and $x$ be a point in the complement component $E_{\alpha}$ of the $V$ amoeba $\mathscr{A}$ in $\mathbb{R}^n$. Let $\Phi_{\alpha}(f)$ be the following expression (see \cite{PR1-04}):
$$
\Phi_{\alpha}(f) = \frac{1}{(2\pi i)^n}\int_{\Log^{-1}(x)}\log
  \left(\frac{f(z)}{z^{\alpha}}\right) \frac{dz_1\wedge \ldots \wedge
    dz_n}{z_1\ldots z_n}.
$$
The function $\log
  \left(\frac{f(z)}{z^{\alpha}}\right)$ has a globally holomorphic branch in the domain $\Log^{-1}(E_{\alpha})$, hence, it defines a holomorphic function in the coefficients of $f$ which takes its values in $\mathbb{C}/ 2\pi i \mathbb{Z}$. In particular, the integral defining $\Phi_{\alpha}(f)$ is independent of the choice of $x$ in $E_{\alpha}$. It has been proved  by Passare and Rullg\aa rd \cite{PR1-04}, that the tropical coefficients of the tropical polynomial defining the spine of the amoeba are those given by the real part of $\Phi_{\alpha}(f)$. Moreover, we can remark in the last example, and generally if $n=1$,\,  the values of  $\Arg (\nu_{(2,\, f)}(\alpha ))$ is precisely the imaginary part of $\Phi_{\alpha}(f)$ modulo $2\pi$. 

\noindent Recall that we assume that $\Delta\cap \mathbb{Z}^n$ contains no point of type I. If $\Ima (\Phi_{\alpha}(f))$ denotes the imaginary part of $\Phi_{\alpha}(f)$, then we define $\nu_{(2,\, f)}$ as follow:
\[
\nu_{(2,\, f)}(\alpha ) = \left\{ \begin{array}{ll}
\exp (i\Ima (\Phi_{\alpha}(f))) &\mbox{if\,  $\alpha\in \Ima (\ord_f)$}\\
1& \mbox{ otherwise}
\end{array}
\right. 
\]


\vspace{0.2cm}

Let  $\Delta_{i}\in \tau$, and  $(a_{i}, b_{i})\in \mathbb{R}^n\times\mathbb{R}$ be the vector  such that the hyperplane in $\mathbb{R}^n\times\mathbb{R}$ containing the extended $\Delta_{i}$ is given by the equation $y = <x , a_{i}> + b_{i} $. 
Let $\{ f_u\}_{u\in ]0;\frac{1}{e}]}$ be the family of polynomials defined by:
\begin{eqnarray}
f_u(z)& =& \sum_{\alpha\in\Ima (\ord_{f})}u^{-c_{\alpha }}
\nu_{(2,\, f)}(\alpha )z^{\alpha}\nonumber\\
&& + \sum_{\alpha\in \Delta_{i} \, of \, type\, II}u^{<\alpha ,a_{i}> +b_{i}+1}z^{\alpha}.\nonumber
\end{eqnarray}

\vspace{0.2cm}

\noindent {\bf  Third step: Geometric and analytic interpretation of $\nu_{(2,\, f)}$ for any $n$.} ${}$

\vspace{0.2cm}

Let's look more closely the function $\Phi_{\alpha}(f)$ defined by Passare and Rullg\aa rd. Let $E_{\alpha}$ be a complement component of the complex amoeba $\mathscr{A}$ of order $\alpha$, and $x\in E_{\alpha}$. We can write $f(z)$ as $\frac{z^{\alpha}f(z)}{z^{\alpha}}$ for any $z\in \Log^{-1}(x)$. On one hand we have:
$$
\log (\frac{z^{\alpha}f(z)}{z^{\alpha}}) = \log\mid f(z)\mid + i\Arg (f(z)),
$$
on the other hand we have :
\begin{eqnarray}
\log (\frac{z^{\alpha}f(z)}{z^{\alpha}})& =&  \log (z^{\alpha}) + \log (\frac{f(z)}{z^{\alpha}}) \nonumber\\
&=& <\alpha ,\, \Log (z) > + i <\alpha ,\, \Arg (z)> + \log (\frac{f(z)}{z^{\alpha}}).\nonumber
\end{eqnarray}
So if we put $\Log (z) = x = (x_1,\ldots ,x_n)$, we obtain:
$$
\mathscr{N}_f(x) \, +\,  i \int_{\Log^{-1}(x)}\Arg (f(z))d\eta (z) = <\alpha ,\, x> + \Phi_{\alpha}(f) + i \int_{\Log^{-1}(x)} <\alpha ,\, \Arg (z)>d\eta (z),
$$
where $d\eta (z) = \frac{1}{(2\pi i)^n}\frac{dz_1\wedge \ldots \wedge
    dz_n}{z_1\ldots z_n}$ and $\mathscr{N}_f$ the Ronkin function. Hence, if we denote by  $m(x)$ the average of $\Arg (z)$  in $\Log^{-1}(x)$ and by  $\mathscr{M}_f(x)$ the average of $\Arg (f(z))$ in $\Log^{-1}(x)$, we obtain:
$$
\mathscr{M}_f(x) = <\alpha ,\, m(x)> + \Ima (\Phi_{\alpha}(f)).
$$
This expression is constant in $E_{\alpha}$ and it depends only on $\alpha$. Recall that $\Ima (\Phi_{\alpha}(f))$ is  precisely the argument of $\nu_{(2,\, f)}(\alpha )$. This determines a collection of hyperplanes $\mathcal{H} = \bigcup_{\alpha\in\Ima (\ord_f)} \mathcal{F}_{\alpha}$ in $\mathbb{R}^n\times\mathbb{R}$ where the first factor (resp. the second factor) is the universal covering of the real torus $(S^1)^n$ where  the coamoeba  sites (resp. the circle $S^1$ where the image of  $\nu_{(2,\, f)}$ sites). The hyperplanes $\mathcal{F}_{\alpha}$ are defined by the   equations $u = \, <\alpha ,\, y> + \Ima (\Phi_{\alpha}(f))$ with $\alpha\in \Ima (\ord_f)$. Another interpretation of $m(x)$ is the following:
consider for each $j=1,\ldots ,n$ a loop $\delta_j$ in the real torus $\Log^{-1}(x)$ on which all the coordinates are constant except the $j^{th}$ coordinate. Let $D_j^x$ be a disc whose boundary is the loop $\delta_j$.

The $j^{th}$ coordinate of of $\alpha$ is precisely the linking number of $V$ with the boundary of the disc $D_j^x$, i.e., the number of intersection points of $V$ with $D_j^x$ counted with signs given by the orientation of $V$. It is Mikhalkin's topological interpretation of Forsberg, Passare, and Tsikh's Theorem 2.1.

\vspace{0.2cm}

Let $z\in V$, $x = \Log (z)$, and   $v_1$ , $v_2$ are points in $\mathbb{R}^n\setminus \mathscr{A}$ such that:
\begin{itemize}
\item[(i)]\, For $i=1,2$, the point $v_i$ is contained in a complement component of order $\alpha_i$;
\item[(ii)]\,  The intersection $[v_1; v_2]\cap \mathscr{A}$ is connected;
\item[(iii)]\, $z$ is contained in the cylinder $\Log^{-1}([v_1; v_2]) \subset (\mathbb{C}^*)^n$ which we denote by 
$(S^1)^n\times [v_1; v_2]$, whose image under the logarithmic map is
precisely the interval $[v_1; v_2]$, and  $\Log ((S^1)^n)\times \{ v_i\}) = v_i$ for $i=1, 2$.
\end{itemize}

Assume that $z=(e^{x_1+i\theta_1},\ldots ,e^{x_n+i\theta_n})$. Let $v_i = (v_{i1},\ldots ,v_{in})$, and $\gamma_{ij}$ be the circle in $\Log^{-1}(v_i)$ parametrized by $(e^{v_{11}+i\theta_1},\ldots , e^{v_{1j}+i\theta_j},\ldots , e^{v_{1n}+i\theta_n})$ with $0\leq \theta_j < 2\pi$. Let $D_j$ be a disc with boundary  $\gamma_{1j}$ containing the cylinder $\gamma_{1j} \times [v_1; v_2]$. Let $\delta_j$ be a circle in the cylinder $\gamma_{1j} \times [v_1; v_2]$ homologous to zero in $(\gamma_{1j} \times [v_1; v_2])\setminus V$, counterclockwise oriented, and the intersection  of its interior with $V$ is empty. Let $T' = \delta_1\times \cdots \times \delta_n$ be the $n$-torus in $(S^1)^n\times [v_1; v_2]$. The cylinder $(S^1)^n\times [v_1; v_2]$ is equipped with the product measure of $d\eta (z)$ on the torus $(S^1)^n$ with the Lebesgue measure on the interval $[v_1; v_2]$; and denote by $d\eta ' (z)$ the restriction of this product  to $T'$.

\vspace{0.2cm}

Let $\varepsilon$ be a small positive real number, and  let us denote by $\delta_{j,\, \varepsilon}$ a family of circles in the cylinder $\gamma_{1j} \times [v_1; v_2]$ which converge  to the union of  $V\cap D_j $ and the two circles $\gamma_{1j}$ and $\gamma_{2j}$ when $\varepsilon$ tends to zero after concatenation. Let ${T'}_{\varepsilon}$  be the torus $\delta_{1,\, \varepsilon}\times\cdots \delta_{n,\, \varepsilon}$. For any $\varepsilon$, the torus  ${T'}_{\varepsilon}$ is  homologous to zero in $(\mathbb{C}^*)^n\setminus V$. Hence, the integrals  $\int_{T'}\Arg (f(z))d\eta ' (z)$ and $\int_{{T'}_{\varepsilon}}\Arg (f(z))d\eta ' (z)$  are equal. So, we have:
$$
\int_{T'}\Arg (f(z))d\eta ' (z) \,  = \,\, \lim_{\varepsilon\rightarrow 0}\int_{T'_{\varepsilon}}\Arg (f(z))d\eta ' (z).  \quad\quad (4)
$$

\vspace{0.1cm}

\begin{Def} Let us denote by $\mathscr{T}_{v_1v_2}$ the set of $n$-torus $T'$ in $((S^1)^n  \times [v_1; v_2])\setminus V$ homologous to zero and constructed as above. The set $Av(v_1,v_2) := \{ \int_{T'}\Arg (z)d\eta ' (z) \in \mathbb{R}^n\, \mid \, T'\in\mathscr{T}_{v_1v_2}\}$ is called the {\em average} of $\Arg (z)$ for $z\in ((S^1)^n  \times [v_1; v_2])\setminus V$.
\end{Def}

We denote by $T_{v_1v_2} := \Log^{-1}(v_1) \cup \Log^{-1}(v_1) \cup (V\cap ((S^1)^n  \times [v_1; v_2]) )$. We say that the sequence $\{T_n'\} \subset \mathscr{T}_{v_1v_2}$ converge to $T_{v_1v_2}$ after concatenation when $n$ tends to infinity if for any $\varepsilon > 0$ there exists $N_0\in \mathbb{N}$ such that if $n>N_0$ and   $z\in T_{v_1v_2}$ then  $\inf_{y\in T_n'}d(z,y)<\varepsilon$, and we denote by 
$\{T_n'\}\stackrel{concat}{\rightarrow}T_{v_1v_2}$.
The set of sequences $\{T_n'\} \subset \mathscr{T}_{v_1v_2}$  such that there exists $j$ and $z\in V\cap (\gamma_j\times [v_1; v_2])$ with $\lim_{n\rightarrow \infty}\inf_{y\in T_n'}d(z,y) = 0$ is denoted by $\mathscr{S}(v_1,v_2)$. We can see easily that the boundary of $Av(v_1,v_2)$ is precisely the set $\{ \lim_{n\rightarrow \infty} \int_{T_n'}\Arg (z)d\eta ' (z) \in \mathbb{R}^n\, \mid \, \{ T_n'\}\in  \mathscr{S}(v_1,v_2)\}$.

\vspace{0.1cm}

\begin{Lem} The boundary of $Av(v_1,v_2)$ is the set 
$$
\mathscr{B}(v_1,v_2) = \{ \lim_{n\rightarrow \infty} \int_{T_n'}\Arg (z)d\eta ' (z) \in \mathbb{R}^n\, \mid \, \{ T_n'\}\in \mathscr{T}_{v_1v_2} \,\, and \,\,
\{T_n'\}\stackrel{concat}{\rightarrow}T_{v_1v_2}  \} .
$$
\end{Lem}

\begin{prooof} The inclusion of $\mathscr{B}(v_1,v_2)$ in the boundary of $Av(v_1,v_2)$ is obvious. Let $a\in \partial Av(v_1,v_2)$, so, $a = \lim_{n\rightarrow \infty} \int_{T_n'}\Arg (z)d\eta ' (z)$ with $\{ T_n'\}\in  \mathscr{S}(v_1,v_2)$. Hence there exist $j$ and $z\in (V\cap (\gamma_j\times [v_1;v_2]))$ with $\lim_{n\rightarrow \infty} d(z,T_n') = 0$. Hence, after concatenation, the projection  on the $j^{th}$ coordinate of the limit of $T_n'$, when $n$ tends to infinity is a loop homologous to  a loop $\delta_j'$ containing the point $z$. We denote by $T_{lim}'$ the limit torus. By hypothesis the integral of $\Arg (z)$ over $T_{lim}'$ exists, it's equal to $a$. But there exists a family  $\{ l_t\}_{t=0}^{t=1}$ of  homologous torus to $T_{lim}'$ such that, $l_0 = T_{lim}'$ and $\lim_{t\rightarrow 1}l_t = T_{v_1v_2} $. So, $a = \lim_{t\rightarrow 1} \int_{l_t}\Arg (z)d\eta ' (z)$, and then $a\in \mathscr{B}(v_1,v_2)$ because the family of the torus is such that $\{l_t\}\stackrel{concat}{\rightarrow}T_{v_1v_2}$.
\end{prooof}

\vspace{0.2cm}

\noindent If $\{ T_n'\}$ is a sequence in $\mathscr{B}(v_1,v_2)$, then  the  limit of the integral $\int_{{T_n'}}\Arg (z)d\eta ' (z)$ when $n$ tends to infinity is denoted by $y$. By Lemma 5.7, those are the  points in the boundary of $Av(v_1,v_2)$. 

\begin{Lem} With the above notations, the set of points in the boundary of  $Av(v_1,v_2)$
 is given by:
$$
\int_{T'}\Arg (f(z))d\eta ' (z) \,  = \,\,   <\alpha_1 - \alpha_2\, , \, y> + \Ima (\Phi_{\alpha_1}(f)) - \Ima (\Phi_{\alpha_2}(f)) + \pi + 2k\pi ,
$$
where $k= \sum_{j=1}^n(lk(V,\gamma_{1j}) - lk(V,\gamma_{2j}))$, and $lk(V,\gamma_{ij})$  denotes the linking number between $V$ and $\gamma_{ij}$ for $i = 1,\, 2$.
\end{Lem}

\begin{prooof}
Indeed,  
the right side of $(4)$ is precisely  $<\alpha_1 - \alpha_2\, , \, y> + \Ima (\Phi_{\alpha_1}(f)) - \Ima (\Phi_{\alpha_2}(f)) + \pi $.
Indeed, when $n$ tends to infinity, the loops $\delta_{1,\, n}$ tends by hypothesis to the two circles $\gamma_{1j}$,\, $\gamma_{2j}$ union the intersection $V\cap D_j$, where  the first circle $\gamma_{1j}$ is clockwise oriented (see Figure 6). So, the torus $T_n'$ converge, after concatenation, to the union of the two boundary torus where the integration on the first torus, is given using the inverse orientation of the loops on each coordinates, and on the second torus, the integration is given using the right orientation of the loops on each coordinates.
\begin{figure}[h!]
\begin{center}
\includegraphics[angle=0,width=0.6\textwidth]{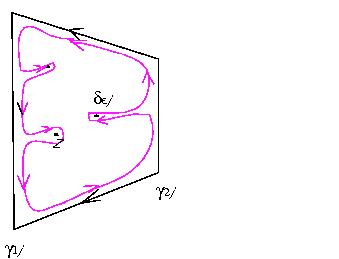}
\caption{}
\label{}
\end{center}
\end{figure}
\end{prooof}

\noindent  In Lemma 5.8, the torus $T'$ is homologous to zero in $((S^1)^n\times [v_1; v_2])\setminus V$. Hence, the average $\int_{T'}\Arg (f(z))d\eta ' (z) $ vanishes, and then, the limit $y$ when $n$ tends to infinity satisfies the following equation:
$$
 <\alpha_1 - \alpha_2\, , \, y> + \Ima (\Phi_{\alpha_1}(f)) - \Ima (\Phi_{\alpha_2}(f)) + \pi  \,\,  =\,\,  0   \mod (2\pi ),\quad\quad (*)
$$

This gives the boundary  of the average of $\Arg (z)$ for $z$ in 
$((S^1)^n  \times [v_1; v_2])\setminus V$; which  is, in the same time,  the boundary of the  average of $\Arg (z)$  for $z$ in  $V\cap ((S^1)^n  \times [v_1; v_2]) $.

\vspace{0.1cm}

Let $\sigma$ be an $(n-1)$-cell in of the spine $\Gamma$ of the amoeba $\mathscr{A}$, and $x\in \sigma$. Let $\mathscr{L}_{(\sigma ,\,  x)}$ be the set of all  lines in $\mathbb{R}^n$ parametrized by a subset $\mathscr{S}_{\sigma}$ of the projective space $\mathbb{P}^1(\mathbb{R}^n)$ satisfying the following:
\begin{itemize}
 \item[(i)]\, Any $l_{(\sigma ,\, x;\, p)} \in \mathscr{L}_{(\sigma ,\,  x)}$ contains the point $x$, and $l_{(\sigma ,\, x;\, p)} $ is not contained in the hyperplane supporting $\sigma$; 
\item[(ii)]\, the intersection of any $l_{(\sigma ,\, x;\, p)} \in \mathscr{L}_{(\sigma ,\,  x)}$ with the complement components $E_\alpha$  and  $E_\beta$  of the amoeba of order respectively $\alpha$ and $\beta$ the ends of the edge $E_{\alpha\beta}$ of $\tau$ dual to $\sigma$ is nonempty.
\end{itemize}

\vspace{0.1cm}

For any $l_{(\sigma ,\, x;\, p)}  \in \mathscr{L}_{(\sigma ,\,  x)} $ , let $v_\alpha\in E_\alpha$ and $v_\beta\in E_\beta$ inside $l_{(\sigma ,\, x;\, p)}$, and $[v_\alpha ;v_\beta ]_p$ be the segment joining $v_\alpha$ and $v_\beta$ passing through $x$. The cylinder $\Log^{-1}([v_\alpha ;v_\beta ]_p) = (S^1)^n\times [v_\alpha ;v_\beta ]_p $ is denoted by $Cyl_{(v_{\alpha} ,v_{\beta} ;x)}^{(p)}$. Let $Cyl_{(\sigma ; x)} := \cup_{p\in \mathscr{S_{\sigma}}} Cyl_{(v_{\alpha} ,v_{\beta} ;x)}^{(p)}$, and $Cyl_{(\sigma)} := \cup_{x\in \sigma} Cyl_{(\sigma ; x)}$, then we have the following:

\begin{Lem} With the above notations, we have:
\begin{itemize}
\item[(i)] \,
$$(\mathbb{C}^*)^n = \bigcup_{\sigma\in \{ (n-1)-cell\,\, of \,\, \Gamma\}} Cyl_{(\sigma)};
$$
\item[(ii)] \, The hyperplane defined by the equation $(*)$ is the same for any cylinder $Cyl_{(v_{\alpha} ,v_{\beta} ;x)}^{(p)}$ with $x\in\sigma$ and $p\in \mathscr{S}_{\sigma}$; i.e., that hyperplane depends only on $\sigma$. On other words, it depends only on $\alpha$ and $\beta$.
\end{itemize}
\end{Lem}

\begin{prooof} The first statement is obvious.
The second statement of the Lemma comes from the fact that the expression of Lemma  5.6 depends only on  the order of the complement components containing the ends of the segment  $[v_{\alpha_1}; v_{\alpha_2}]$ and not on the segment itself.
\end{prooof}

\vspace{0.1cm}

\begin{Def} We denote by $\mathscr{H}$ the union of the hyperplanes defined by the equations $(*)$, and we call it the {\em shell}  of the coamoeba of $V$ (or the {\em average contour} of the coamoeba).
\end{Def}

\vspace{0.1cm}

\begin{Rem} ${}$ The hyperplanes defined in $(*)$ are the codual hyperplanes to the set of edges of $\tau$.
\end{Rem}

\vspace{0.2cm}

\begin{Def}
Let $H_{\alpha\beta}$ be a hyperplanes given by $(*)$, if we know the position of the average of $\Arg (z)$ for $z\in ((S^1)^n\times [v_{\alpha} ;\, v_{\beta}])\setminus V$ relatively to $H_{\alpha\beta}$, then we call the points in the other side of $H_{\alpha\beta}$ {\em the average of the coamoeba} of $V$.
\end{Def}

For any  vertex $v$ of the spine $\Gamma$, we have a subset $\mathscr{H}_v$ of $\mathscr{H}$, which is the union of the hyperplanes given by $(*)$
 corresponding to the edges of $\tau$ dual to the $(n-1)$-cells of $\Gamma$ adjacent to $v$. The side of  the position of the average of $\Arg (z)$ for $z\in ((S^1)^n\times [v_{\alpha} ;\, v_{\beta}])\setminus V$, with $E_{\alpha\beta}$ an external edge, relatively to the hyperplane $H_{\alpha\beta}$, is well defined by Lemma 5.13; and then the average of the coamoeba is also well defined relatively to $H_{\alpha\beta}$.

\vspace{0.2cm}

\noindent {\bf  Fourth step: Determination of the coamoeba average relatively to $\mathscr{H}$.} ${}$

\vspace{0.2cm}

\noindent The purpose of this step is the determination of the position of the average of the coamoeba of $V$ relatively to its shell $\mathscr{H}$. First of all, we determine the position of the average relatively to the external hyperplanes i.e., the hyperplanes codual to the external edges of the Newton polytope. We denote by $E_{\alpha}$ the connected component of the complement of the amoeba of order $\alpha$. We can remark that the position of the average of the coamoeba relatively to  hyperplanes in  $\mathscr{H}$ is then well defined by that of its position relatively to the external codual hyperplanes. We have the 
  the following  Lemma:

\begin{Lem} Let $K$ be a compact in $\mathbb{R}^n$ containing the vertices of the spine $\Gamma$  of $\mathscr{A}$, and $C$ be an $(n-1)$-cell of $\Gamma$ dual to the edge $E_{\alpha\beta}$. Let $U_{\varepsilon}(C)$ be an $\varepsilon$-neighborhood of $(\mathbb{R}^n\setminus ( K\cap \mathscr{A})) \cap  C$.
Let $H_{\alpha\beta}$ be a hyperplane  in the shell of the coamoeba and codual to $E_{\alpha\beta}$.
 Then the position of  the average of $\Arg (\Log^{-1}(U_{\varepsilon}(C))\cap V)$ denoted by $\mathscr{U}(C,V)$,
relatively to $H_{\alpha\beta}$ cannot be situated in the two  side of this hyperplane in the following sense: there exists a subset $S$ in $H_{\alpha\beta}$ with vanishing measure such that for any $x\in H_{\alpha\beta}\setminus S$ there exists an open set $\mathcal{V}(x)$ in $\mathbb{R}^n$ containing $x$ such that $\mathcal{V}(x)\cap \mathscr{U}(C,V)$ is in one side of $H_{\alpha\beta}$.
\end{Lem}

\begin{prooof} 
The set of arguments of points in $V$ with infinite module, looks like the arguments of infinite points in a complex tropical hypersurface. Indeed, if $z\in V$ is such that $|z|$ is so large and $\Log (z) \in U_{\varepsilon}(C)$, then the dominating monomials of the polynomial are those of indices $\alpha$ and $\beta$. But in this case, we can apply the tropical localization, which brings us back to the coamoeba of complex tropical hypersurface; and the non-Archimedean coamoebas is obviously in one side of $H_{\alpha\beta}$ for $z\in \Log^{-1}(U_{\varepsilon}(C))\cap V$. 
\end{prooof}

\begin{Cor} The average of the coamoeba of $V$ is equal to the coamoeba of the complex tropical hypersurface $\lim_{u\rightarrow 0}H_u(V_{f_u})$.
\end{Cor}

\begin{prooof} For any $\mathscr{H}_v$ with $v$ dual to an element $\Delta_i$ of $\tau$ of maximal dimension and having an external edge, the position of the average of the coamoeba of $V$ is well defined by Lemma 5.13. So, it is well defined relatively to  $\mathscr{H}_w$ for $w$ a vertex of $\Gamma$ dual to an element of $\tau$ with common facet with $\Delta_i$. The triangulation is finite, so the position of the average of the coamoeba of $V$ is well defined relatively to all $\mathscr{H}$. By construction of the complex tropical  hypersurface $\lim_{u\rightarrow 0} H_u(V_{f_u})$, the codual hyperplanes arrangement are the same of that given by the equations $(*)$, because these hyperlanes are defined by the same equations modulo $2\pi$. Hence, the coamoeba of that complex tropical  hypersurface is equal to the average of the coamoeba of $V$ defined in 5.12.
\end{prooof}

\vspace{0.1cm}

Corollary 5.14  is the analogous to Passare and Rulld\aa rd Proposition 1  \cite{PR1-04} on amoebas, of course if we see to  the spine as an average of the amoeba. 

\vspace{0.1cm}

\begin{Lem} With the above notations, 
  the closure in the real torus of the coamoeba of $V$ is homeomorphic to the coamoeba of the complex tropical hypersurface $V_{\infty ,\, f} := \lim_{t\rightarrow 0}H_t(V_{f_t})$.
  \end{Lem}
\begin{prooof}
We have to prove that the cardinality of the connected components of the complement of the coamoeba $(S^1)^n\setminus \overline{co\mathscr{A}}$
 is equal to that of the complement of its  average. Let $\mathscr{E}^c$ be a complement component of $\overline{co\mathscr{A}}$, then it is clear that the average coamoeba cannot cover all $\mathscr{E}^c$. On the contrary, this means that the codual hyperplanes containing the subset of critical values of the argument map in the boundary of $\mathscr{E}^c$, intersect in one point, and then $\mathscr{E}^c$ is empty. Contradiction,  because if it is nonempty, then  its boundary must contains more than one point in the critical values of the argument map and in the intersections of the codual hyperplanes. So the cardinality of the connected components of the complement of the coamoeba $(S^1)^n\setminus \overline{co\mathscr{A}}$ cannot exceed that of the average coamoeba. Let $\mathscr{E}^c_{av}$ be a complement component of the average coamoeba, its boundary contains a subset of the boundary of some complement component of the complex coamoeba (the points in the intersections of the codual hyperplanes). Hence, the cardinality of the connected components of the complement of the  average coamoeba $(S^1)^n\setminus co\mathscr{A}_{av}$ cannot exceed that of the complex coamoeba. The Lemma is done, because  the average of the complex coamoeba  is precisely the coamoeba of the complex tropical hypersurface $V_{\infty ,\, f} := \lim_{t\rightarrow 0}H_t(V_{f_t})$. 

\end{prooof}

\vspace{0.1cm}

{\bf Geometric interpretation of indices of type I.}${}$

\vspace{0.2cm}

We give in this paragraph, the interpretation of  type I indices in terms of normal vectors to the facets of the extended Newton polytope $\Delta$  corresponding to the elements of maximal dimension of the subdivision $\tau$ of $\Delta$ adjacent that indices after a small perturbation of the corresponding coefficient. This gives a geometric interpretation of the valuation of those indices. On the other hand we give a geometric interpretation of their argument part, in terms of the codual hyperplanes sub-arrangement to the adjacent edges to these indices.

\vspace{0.2cm}

\noindent {\it End  of Theorem 5.2 proof.} Let us define $\nu_{(2,\, f)} (\beta )$ for $\beta$ of type I. We know that the function $\Phi_{\beta}$ viewed as function on the coefficients of the polynomial, is holomorphic function (see \cite{PR1-04}, and \cite{R-01}). So, if we replace the coefficient $a_{\beta}$ of $f$ by  $(1-\varepsilon )a_{\beta}$ such that the new polynomial $f_{\varepsilon}$ has $\beta$ in the image of its order map, then the limit of $\Phi_{\beta}(f_{\varepsilon})$ when $\varepsilon$ tends to zero is well defined.  Then, we put $\Phi_{\beta}(f) $ equal to that limit. In this case,  Lemma 5.15 is still valid, and   Theorem 5.2  is an immediate consequence of  Lemma 5.15.

\noindent  Recall that Passare and Rullg\aa rd \cite{PR1-04} prove the analogous for the spine of the amoeba; in this case, the image of  $\alpha\in\Ima (\ord )$ by the corresponding function, which is $\nu_{PR}$ is given by  $\R (\Phi_{\alpha}(f))$. Moreover, they prove that if $\alpha \in \Verte (\Delta )$ then $\R (\Phi_{\alpha}(f)) = \log\mid a_{\alpha}\mid$.

\vspace{0.2cm}

\begin{Def} Let $V$ be a complex algebraic hypersurface defined by a polynomial $f$. Let $\mathscr{H}_f$ be its associated codual hyperplanes arrangement, and $\Critv (\Arg)$ be the set of critical values of the argument map. Hence, the set $\overline{co\mathscr{A}_f}\setminus \left( \Critv (\Arg) \cup\, \mathscr{H}_f\right) = \bigsqcup_{j=1}^l \mathcal{R}_j $. The {\em extra-pieces} are the $\mathcal{R}_j$'s such that their boundary is not contained in the arrangement $\mathscr{H}_f$. 
\end{Def}

\vspace{0.2cm}

\begin{Rems}${}$
\begin{itemize}
\item[(1)]\, It is easy to see that for given polytope $\Delta$, if the number of complement components of the closure in the real torus $(S^1)^n$ of a  complex tropical coamoeba is maximal, then it can be realized by a maximally sparse polynomial such that the spine of its amoeba has one vertex. In this case, the boundary of any complement component of the coamoeba is contained in the codual hyperplane arrangement $\mathscr{H}$ associated to the trivial subdivision of $\Delta$, this means that  $\mathscr{H}$ is precisely the union of the codual hyperplanes to the edges of $\Delta$.
\item[(2)]\,  If $n=2$, then the area  (counted with multiplicity) of complex algebraic plane  curve coamoebas with Newton polygon $\Delta$, cannot exceed $2\pi^2\Area (\Delta )$. Indeed, for the complex tropical plane curves it is clear.  Let $v$ be a vertex of the non-Archimedian amoeba of that complex tropical curve $V_{\infty}$; the deformation of $\Arg (\Log^{-1}(v)\cap\, V_{\infty})$ given in Theorem 5.2  gives rise to a deformation of the codual hyperplanes arrangement $\mathscr{H}_v$, such that when $t=\frac{1}{e}$
  the arrangement $\mathscr{H}_{v,\, \frac{1}{e}}$ has the same combinatorial type of $\mathscr{H}_v$ (by construction). Hence, the set of arguments limited by $\mathscr{H}_v$ is deformed to a set $S$ limited by $\mathscr{H}_{v,\, \frac{1}{e}}$ maybe with some extra-pieces. The area of $S$ cannot exceed $2\pi^2\Area (\Delta_v)$ (here $\Delta_v$ is the dual, of the vertex  $v$, in the subdivision of the Newton polygon). Indeed, there are two cases: (i) if there is no coefficients in $\Delta_v$ other than of its vertices it's clear, (ii) if there are some coefficients in $\Delta_v$, then we can see, by deformation, that the area of $S$ cannot exceed $2\pi^2\Area (\Delta_v)$. 
It suffice to consider the deformation of a polynomial having those coefficients in the image of its order map, and the deformation is given only on the module of the coefficients with index not in $\Verte (\Delta_v)$ (the argument of the coefficients of such deformation are those of the initial polynomial, such deformation exists, see \cite{R2-00}). More precisely, for any time $t$ of the deformation,  the set of arguments is limited by the arrangement $\mathscr{H}_t$, and then its area cannot exceed that of the coamoeba of the complex tropical hypersurface $V_{\infty ,\, t}$.
\end{itemize}
\end{Rems}

\begin{Rem}
Let $V$ be a complex algebraic hypersurface, then by Theorem 5.2 there exists a complex tropical hypersurface $V_{\infty}$ such that:
\begin{itemize}
\item[(a)]\, The spine $\Gamma$ of the amoeba of $V$ is the same as the non-Archimedean amoeba $\Gamma_{\infty}$ which the image under the logarithmic map of $V_{\infty}$;
\item[(b)]\, The closure in the real torus of the amoebas of $V$ and $V_{\infty}$ are homeomorphic;
\item[(c)]\, The volume of the coamoeba of $V_{\infty}$ (counted with multiplicity) is greater or equal to that of the coamoeba of $V$, and equality hold if and only if the set of critical values of the argument map restricted to $V$ is discrete, which means that its coamoeba contains no extra-piece;
\item[(d)]\, If $n>1$, then the codual hyperplanes arrangement associated to $V_{\infty}$ determines completely the topology of the coamoeba of $V$, but not its geometry (because of the extra-pieces which can exist in the classical complex case).
\end{itemize}
\end{Rem}

\newpage

\begin{The}
Let $V$ be an algebraic hypersurface in $(\mathbb{C}^*)^n$ defined
by a polynomial $f$ with Newton polytope $\Delta$ and coamoeba $co\mathscr{A}$.
Then we have:
\begin{itemize}
\item[(a)]\, The interior of any connected component of $(S^1)^n\setminus  co\mathscr{A}$  is a convex set (if nonempty);
\item[(b)]\, The number of connected components of $(S^1)^n\setminus \overline{  co\mathscr{A}}$
is not greater than $n! \Vol (\Delta )$ where $\overline{  co\mathscr{A}}$ is the closure of
$co\mathscr{A}$ in the flat torus  $(S^1)^n$.
\end{itemize}
\end{The}

\begin{prooof}
Using the analysis tools, more precisely  Bochner Theorem, M. Passare  has proved the convexity  of the components of the coamoebas complement (\cite{P1-07}, private communication).\\
The proof of the second statement of this Theorem use firstly the fact that the closure of the complex coamoeba in the real  torus  is homeomorphic to that of a complex tropical coamoeba, secondly    we use the  remark 5.5 (1), induction on the volume of $\Delta$, and the convexity of $\Delta$. Firstly, if $\Delta$ is a simplex, then the number of connected components of $(S^1)^n\setminus \overline{co\mathscr{A}}$ is  $n! \Vol (\Delta )$. Let $\Delta = \Delta_1  \cup \, \Delta_2$
 such that $\Delta_1  \cap \, \Delta_2$ is a face of dimension $n-1$, and $\Delta_j$ are with nonempty interior for $j=1,\, 2$. Using induction on the volume, the number of connected components of $(S^1)^n\setminus \overline{  co\mathscr{A}_{\Delta_j}}$ is not greater than $n! \Vol (\Delta_j )$ for $j= 1 , 2$. We start by  the following definition and Lemma:

\vspace{0.1cm}

\begin{Def}A point $p$ in the coamoeba  of some complex tropical hypersurface is called {\em separating point} if $p$ is contained in the closure of two different complement components of its coamoeba.
\end{Def}

By remark (2) in  4.1 we know that the position of the coamoeba relatively to the external codual hyperplanes is well determined by their  canonical orientation given as described in that remark. 

\begin{Lem} Let $\mathscr{E}^{(c,\, j)}_{\Delta_1}$ be  two connected components of the complement of the coamoeba with degree $\Delta_1$ ,\,  $j=1,\, 2$, and  $\mathcal{B}_j$  are the boundary of $\mathscr{E}^{(c,\, j)}_{\Delta_1}$ such that $\mathcal{B}_1\cap\, \mathcal{B}_2 = \{ p \}$ with  $p$ a separating point.  Then at most one of the $\mathcal{B}_j$'s can intersect  the same two  different complement components of the coamoeba of degree $\Delta_2$
\end{Lem}

\begin{prooof} If there exist  $\mathcal{B}_1$ and  $\mathcal{B}_2$ satisfying the hypothesis of the Lemma, then there exist at least two edges $E_1$ of $\Delta_1$ and $E_2$ of $\Delta_2$ such that the convex hull of $\partial E_1 \cup\,\partial E_1 $ is not contained in $\Delta_1\cup\,\Delta_2$. Indeed, using the natural orientation of the external codual hyperplanes corresponding to $\mathcal{B}_1$ and  $\mathcal{B}_2$ we obtain the position of the coamoeba relatively to those hyperplanes. The fact that they intersect the same two different complement components of the coamoeba of degree $\Delta_2$, imply that there is another $\mathcal{B'}_1$ and  $\mathcal{B'}_2$ which are the boundary of two connected components of the complement of the coamoeba with degree $\Delta_2$ which we denote by  $\mathscr{E'}^{(c,\, j)}_{\Delta_2}$, and satisfying the same hypothesis of the Lemma (see Figure 7 and 8 for $n=2$). This means that $\mathcal{B}_1$,\, $\mathcal{B}_2$,\,  $\mathcal{B'}_1$, \, and \, $\mathcal{B'}_2$  are codual to some  edges such that the convex hull of there union is not contained in the union of $\Delta_1$ and $\Delta_2$.

\begin{figure}[h]
\includegraphics[angle=0,width=0.2\textwidth]{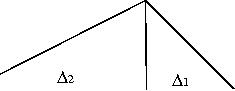}\quad\quad
\includegraphics[angle=0,width=0.3\textwidth]{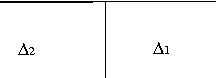}\quad\quad
\includegraphics[angle=0,width=0.3\textwidth]{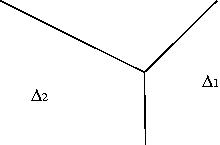}
\caption{The left picture represent the case when $\Delta_1\cup\, \Delta_2$ is convex, on the middle the case when the union is convex and they have two facets where each pair is contained in the same hyperplane, on the right the case where the union is no-convex.}
\label{}
\end{figure}

\begin{figure}[h]
\includegraphics[angle=0,width=0.2\textwidth]{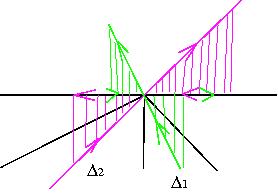}\quad\quad
\includegraphics[angle=0,width=0.3\textwidth]{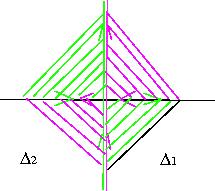}\quad\quad
\includegraphics[angle=0,width=0.3\textwidth]{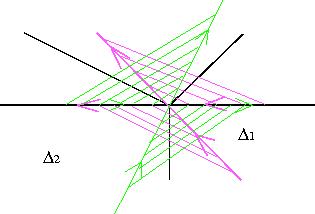}
\caption{The position of the  coamoeba in the three cases illustrated in figure 6.}
\label{}
\end{figure}

\end{prooof}

\begin{figure}[h]
\includegraphics[angle=0,width=0.2\textwidth]{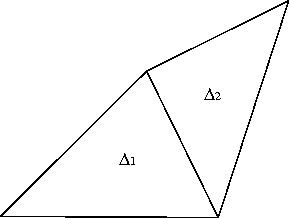}\quad\quad\quad
\includegraphics[angle=0,width=0.3\textwidth]{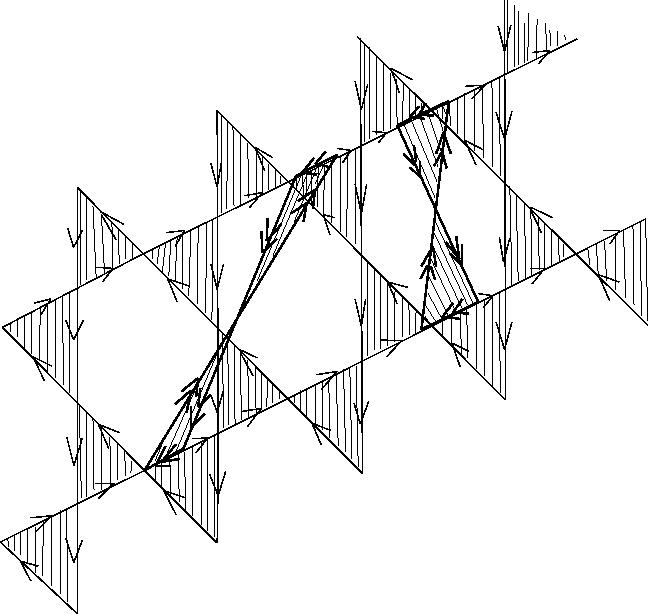}
\caption{In the right figure, only the right double triangles can occur in the convex case (i.e., the union
$\Delta_1\cup\Delta_2$ is convex), and the left double triangles can be in the non convex case.}
\label{}
\end{figure}

{\it End of the proof of statement (b) of Theorem 5.19.} By Lemma 5.21, the number of complement component of the coamoeba of degree $\Delta$ is at most the number of complement component of the coamoeba of degree $\Delta_1$ plus the number of complement component of the coamoeba of degree $\Delta_1$ . Indeed, $\mathcal{B}_1$ and  $\mathcal{B}_2$ satisfying the hypothesis of Lemma 5.21 cannot separate the complement components of the coamoeba of degree $\Delta_2$ to more than $n!\Vol (\Delta_2) + 1$. 
\end{prooof}

 The first statement (a) of Theorem 5.19,  results from the next Lemma 5.22.
First of all, we know that the convexity is a local property, so we have to prove the first part of Theorem 5.19 only locally. Moreover, our main result deals with the two dimensional case, but the argument used works in all dimension. This argument is strongly inspired to Mikhalkin's works in \cite{M3-00}, and more precisely Lemma 1 of that paper. Therefore, the first statement of our Theorem 5.19 is an immediate consequence of the following Lemma: 

\begin{Lem} Let $z\in V_f$ be a critical point of  $\Arg_{\mid V_f}$ the restriction of the argument map to $V_f$, $D$ be a disk in the torus $(S^1)^2$ containing $\Arg (z)$, and $B$ be a component of $(\Arg_{\mid V_f})^{-1}(D)$ containing $z$. Then $D\setminus \Arg (B)$ is convex.
\end{Lem}

\begin{prooof}
Assume on  the contrary that there exists a component of the coamoeba complement which is not convex. Hence there exists a closed interval $I\subset D$ such that its both endpoints are in the same component of $D\setminus \Arg (B)$ such that $I$ has a rational slope, and the intersection of $\Arg^{-1}(I)$ with $B$ is not empty (see figure 7).  The fact that $D$ is contractible imply that the intersection $B\cap \Arg^{-1}(I) $ contains no cycle in $B$, this means that it has a topology of a segment. On the other hand, the slope of $I$ is rational, so we can choose a non vanishing complex numbers $a$ and $b$ such that the holomorphic annulus $\mathscr{H} = \{ z\in (\mathbb{C}^*)^2 \mid az^{\alpha}+bz^{\beta} = 0 \}$ contains $\Arg^{-1}(I)$, project properly to $I$ and intersect $B$ ($\alpha$ and $\beta$ are integer numbers well determined by the slope of $I$). Hence we have a contradiction because $B\cap \Arg^{-1}(I)$ contains no cycle in $(\mathbb{C}^*)^2$, which means that $B\cap \Arg^{-1}(I)$ is empty. Because $B$ and $\Arg^{-1}(I)$ are holomorphic and then their intersection (if nonempty) contains a cycle.
\end{prooof}

\begin{figure}[h!]
\begin{center}
\includegraphics[angle=0,width=0.4\textwidth]{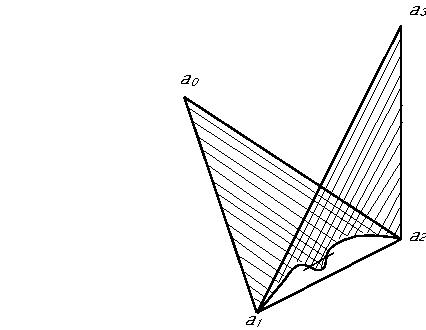}
\caption{}
\label{c}
\end{center}
\end{figure}

\begin{Rem}
Using  the fact that the sets of critical points of the logarithmic map and the argument map are the same, and combining it with the proof of Lemma 1 in  \cite{M3-00}, we obtain the result.
\end{Rem}

\vspace{0.2cm}

\section{Characterization of the virtual vanishing coefficients}

\vspace{0.2cm}

Let $f$ be a polynomial as in $(1)$, and denote by $A := \supp (f) \cup\, \Ima (\ord_f)$, and by $\tau$ the subdivision of the Newton polytope $\Delta$ (which we assume a triangulation) of $f$ dual to the spine of the amoeba $\mathscr{A}$ of the complex hypersurface $V$ with defining polynomial $f$. Also we denote by $\nu_{(2,\, f)}$ the function defining the codual hyperplanes arrangement $\mathscr{H}_f$ associated to the coamoeba $co\mathscr{A}$ of $V$. Let $\nu : A \rightarrow \mathbb{R} \times S^1$ be the function defined as follow:
\begin{itemize}
\item[(i)]\, if $\alpha\in \Ima (\ord_f)$, then  $\nu (\alpha
  )= (-c_{\alpha},\nu_{(2,\, f)}(\alpha))$
\item[(ii)]\, Let
  $\Delta_i$ be an element of $\tau$ with maximal
  dimension, and  $y=<x,a_i>+b_i$ be the equation of the hyperplane in
  $\mathbb{R}^n\times \mathbb{R}$ containing the points of coordinates
  $(\alpha , -c_{\alpha})\in\mathbb{R}^n\times \mathbb{R}$ for
  $\alpha\in \Verte (\Delta_i )$, \, $a_i=(a_{1,\, i},\ldots ,a_{
  n,\, i})\in \mathbb{R}^n$ and $b_i\in\mathbb{R}$. 
If $\alpha\in\Delta_i\setminus \Ima (\ord_f)$, and for any small perturbation of the coefficient $a_{\alpha}$ the amoeba of the new polynomial remains without complement component of order $\alpha$, then   
we set $\nu (\alpha ) = (<\alpha ,a_i>+b_i +1, 1)$, this means $\alpha$ is of type II. Otherwise we set $\nu (\alpha ) = (<\alpha ,a_i>+b_i,\nu_{(2,\, f)}(\alpha ))$. We denote by $\nu_1$ the first coordinate of $\nu$ 
\end{itemize}

In this section $f_t$ is the family of polynomials  defined by:

$$
f_t(z)= \sum_{\alpha\in \supp (f)} a_{\alpha}t^{\nu_1(\alpha )}e^{\nu_1(\alpha )}z^{\alpha},\qquad\qquad\qquad\qquad\qquad   (5)
$$
and we denote by $a_{(\alpha ,\, t)}$ the coefficient of $f_t$ with index $\alpha$.

\vspace{0.1cm}

\noindent Let 
$\{ \tilde{f}_u\}_{u\in ]0,\, \frac{1}{e}]}$ be the family of polynomial defined by:
$$
\tilde{f}_u(z) = \sum_{\alpha\in A}e^{\nu_1 (\alpha )}u^{\nu_1 (\alpha )}(eua_{\alpha} + (1-eu)\nu_{(2,\, f)}(\alpha ))z^{\alpha}.\qquad\qquad\qquad\qquad   (6)
$$
By Theorem 5.2,  the coamoeba of the hypersurface defined by $f$ and the coamoeba of the complex tropical hypersurface $V_{\infty ,\, f} := \lim_{u\rightarrow 0} H_u(V_{\tilde{f}_u})$ are homeomorphic. More generally, let $A_t = \supp (f) \cup\, \Ima (\ord_{f_t})$, and $\nu_1^{(t)}$ (resp. $\nu_{(2,\, f_t)}$ ) be the function defining the spine of the amoeba $\mathscr{A}_t$ of $f_t$ (resp. the arrangement of hyperplanes in the average of the coamoeba of $f_t$). Consider the family of polynomials $\tilde{f}_u^{(t)}$ defined by : 
$$
\tilde{f}_u^{(t)}(z) = \sum_{\alpha\in A_t}e^{\nu_1^{(t)} (\alpha )}u^{\nu_1^{(t)} (\alpha )}(eua_{(\alpha ,\, t)} + (1-eu)\nu_{(2,\, f_t)}(\alpha ))z^{\alpha}.\qquad\qquad\qquad\qquad   (7)
$$
Then for any $t$ the coamoeba of the hypersurface defined by the polynomial $f_t$,
satisfies the same conclusion as the coamoeba of $f$, that it is homeomorphic to the coamoeba of the complex tropical hypersurface
$V_{\infty ,\, f_t} := \lim_{u\rightarrow 0} H_u(V_{\tilde{f}_u^{(t)}})$. We denote by $\mathscr{H}_f$ the shell of the coamoeba of $V$.

\vspace{0.2cm}

\noindent If $\beta$ is a lattice point in $\Ima (\ord_f)\setminus \supp (f)$, then we denote by $\mathscr{H}_{(\beta ,\, f)}$ the sub-arrangement of $\mathscr{H}_f$ consisting of  the hyperplanes  codual to the edges of $\tau$ adjacent to $\beta$. The following theorem gives a combinatorial characterization which must satisfies the set $\Ima (\ord_f) \cap\, \supp (f)$ for the existence of such $\beta$. As an application of this result,
 we prove that the amoeba of maximally sparse polynomial is solid (we can see the first proof in \cite{N1-06}).
 Moreover, it gives a geometric properties satisfied by the coamoeba of hypersurface defined by polynomial $f$ such that  $\Ima (\ord_f) \setminus\,  \supp (f)$ is nonempty. In particular, if $n=2$, we give a characterization on the number of the coefficients of  real polynomials defining a Harnack curves.

\vspace{0.2cm}

\begin{The} Let $\beta\in \Ima (\ord_f)\setminus \supp (f) $. Then we have the following:
\begin{itemize}
\item[(i)]\, There exists an effective sub-arrangement $\mathscr{H}'$ of $\mathscr{H}_f$ with the same combinatorial type of $\mathscr{H}_{\beta}$ such that any hyperplane $H\in \mathscr{H}'$ is of weight at least two.
\item[(ii)]\, The coamoeba $co\mathscr{A}$ of the hypersurface $V$ with defining polynomial $f$ contains some extra-pieces. In particular, if $n=2$, this imply that the real part $\mathbb{R}V$ of $V$ is not a Harnack curve.
\end{itemize}
\end{The}

\begin{prooof} Let us prove the first statement of the Theorem. 
We denote by $\tau^{(t)}_{\infty}$ the subdivision (which we assume a triangulation) of $\Delta$ dual to  the tropical hypersurface $\Gamma_{\infty ,\, f_t} = \Log (V_{\infty ,\, f_t})$, and by $\mathscr{H}^{(t)}_{\infty}$ the arrangement of codual hyperplanes to the edge of $\tau^{(t)}_{\infty}$ associates to $\nu_{(2,\, f_t)}$.  Let $\tau_{f_t}$ be the subdivision of $\Delta$ dual to the spine of the amoeba $\mathscr{A}_{f_t}$ of $V_{f_t}$.
Let $t_{(max,\, \beta)}$ be the real number in $[0; \frac{1}{e}]$ defined as follow:
$$
t_{(max,\, \beta)} := \max \{ 0, \max \{ s\in ]0;\frac{1}{e}] \mid \, \mathscr{H}_{(\beta ,\, f_s)}\,\, is\,\, effective  \}\} .
$$
This means that any hyperplane $H$ in $\mathscr{H}_{(\beta ,\, f_{t_{(max,\, \beta)}})}$ is codual to an effective edge $E$ of $\tau_{t_{(max,\, \beta)}}$, and $E$ is not adjacent to $\beta$, because in our case $\beta \notin \supp (f)$. For any $t_{(max,\, \beta)}< t \leq \frac{1}{e}$, there exists a virtual sub-arrangement $\mathscr{H}_{(\beta ,\, f_t)}$ of the arrangement $\mathscr{H}_{ f_t}$ codual to the edges of $\tau_{f_t}$ adjacent to $\beta$. Indeed, if we fix $t$, in the deformation of $V_{f_t}$ given by the family $\tilde{f}_u^{(t)}$, the spines  (resp. the codual hyperplanes arrangement) of the amoebas of $V_{\tilde{f}_u^{(t)}}$ (resp. coamoebas) for $u$ sufficiently close to zero and $u$ sufficiently close to $\frac{1}{e}$, they   have the same combinatorial type, because the function $\nu_1^{(t)}$ (resp. $\nu_{(2,\, f_t)}$) defining the spine (resp. the codual hyperplanes arrangement) are the same for $u$ sufficiently close to zero and $u$ sufficiently close to $\frac{1}{e}$ by construction. Moreover, the sub-arrangement of hyperplanes codual to the edges adjacent to $\beta$ are effective for $u\ne \frac{1}{e}$ and virtual for $u = \frac{1}{e}$ by construction. The  compactness of the real torus, and the fact that $\mathscr{H}_{(\beta ,\, f_t)}$ is invariant under the translation group $(2\pi \mathbb{Z})^n$,  imply the convergence of the $\mathscr{H}_{(\beta ,\, f_t)}$ to some  effective sub-arrangement $\mathscr{H}_{(\beta ,\, f_{t_{(max,\, \beta)}})}$ when $t$ tends to $t_{(max,\, \beta)}$, by definition of $t_{(max,\, \beta)}$. Geometrically, this means that the complement components of the amoebas of $V_{f_t}$ of order $\beta$ disappear when $t=t_{(max,\, \beta)}$, moreover, the codual hyperplanes sub-arrangement $\mathscr{H}_{(\beta ,\, f_t)}$ converge to some effective sub-arrangement $\mathscr{H}_{(\beta ,\, f_{t_{(max,\, \beta)}})}$. Hence, any hyperplane $H\in \mathscr{H}_{(\beta ,\, f_{t_{(max,\, \beta)}})}$ is codual to an effective edge $E$ of $\tau_{t_{(max,\, \beta)}}$ such that $E$ is not adjacent to $\beta$, because $\beta\notin \supp (f)$.  So, the hyperplane $H$ is counted at least twice, and then its multiplicity is at least equal two. Let $\mathscr{E}$ be the dual of $\mathscr{H}_{(\beta ,\, f_{t_{(max,\, \beta)}})}$ in the set of effective edges of $\tau_{t_{(max,\, \beta)}}$. By duality, if we denote by $C$ the convex hull of the set of vertices of  $\mathscr{E}$, then $C\setminus \Verte (C)$ contains at least an element of the support of $f$. This means that if $f$ is maximally sparse, then $\Ima (\ord_f)\setminus \supp (f)$ is empty.

\noindent The second statement of Theorem 6.1 result from the fact that
with our  hypothesis, 
there exists a collection of effective edges $\mathscr{E}$ in $\tau_f$ dual to $\mathscr{H}_{(\beta ,\, f)}$ such that: there exist two elements $\Delta_1$ and  $\Delta_2$ in $\tau_f$ of maximal dimension,  and $H\in \mathscr{H}_{(\beta ,\, f)} $  satisfying :
\begin{itemize}
\item[(1)]\, $\Delta_1$ has $\beta$ as vertex, and  with $H$ codual to some edge of $\Delta_1$;
\item[(2)]\, $\Delta_2$ has an effective edge in $\mathscr{E}$ dual to $H$;
\item[(3)]\, The set  $\Arg (V_{\infty ,\, f} \cap \Log^{-1}(v_1)) \,\cap \Arg (V_{\infty ,\, f} \cap \Log^{-1}(v_2))$ has nonempty interior, and the boundary of this intersection contains a subset of $H$ with nonempty interior (here we mean the interior of the hyperplane $H$), see Figure 11 on the right for $n=2$. 
\end{itemize}
Hence, in  the  complex hypersurface $V$, the set of critical values of the argument map is not contained in the union of hyperplanes, and then the extra-pieces exist.

\begin{figure}[h!]
\begin{center}
\includegraphics[angle=0,width=0.4\textwidth]{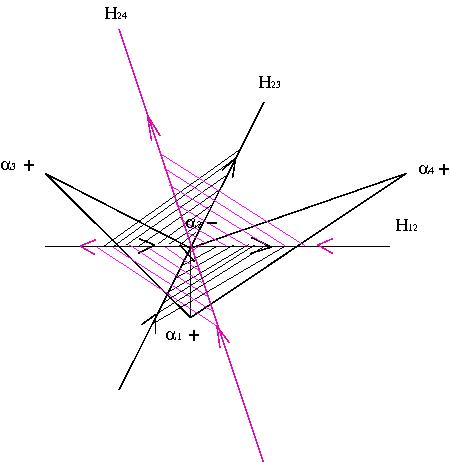}
\includegraphics[angle=0,width=0.5\textwidth]{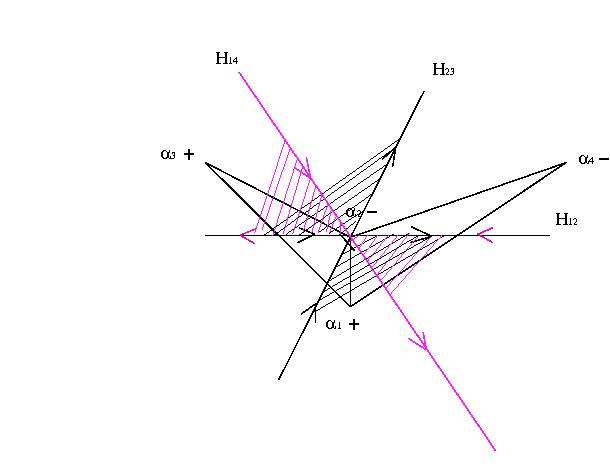}
\caption{The left picture is an illustration of the coamoeba in the case of Harnack curves, and on the right the coamoeba when we must have some extra-pieces.}
\label{}
\end{center}
\end{figure}

\begin{figure}[h!]
\begin{center}
\includegraphics[angle=0,width=0.35\textwidth]{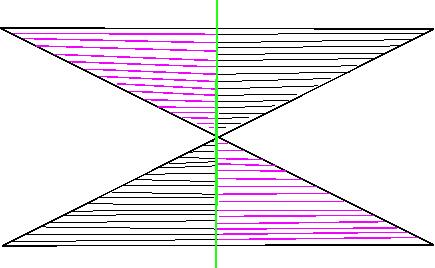}\quad
\includegraphics[angle=0,width=0.35\textwidth]{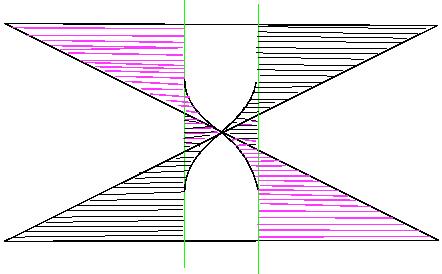}\quad
\includegraphics[angle=0,width=0.2\textwidth]{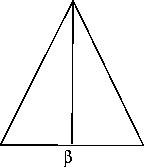}
\caption{Example of the deformation of the  codual hyperplanes  arrangement when $\beta\in F$ with $F$ a facet of $\Delta_i$ with change of combinatorial type.}
\label{}
\end{center}
\end{figure}

\end{prooof}

\vspace{0.1cm}

\begin{The}[Nisse \cite{N1-06}] The amoeba of a complex algebraic hypersurface  defined by  maximally sparse polynomial is solid.
\end{The}

\begin{prooof}
Indeed, if the polynomial $f$ is maximally sparse, then the first statement of Theorem 6.1 implies that the set  $\Ima (\ord_f)\setminus \supp (f)$ is empty.
\end{prooof}

\vspace{0.1cm}

\begin{The} Let $V$ be a complex  algebraic  plane curve defined by a polynomial $f$ with Newton polytope $\Delta$, such that its real part $\mathbb{R}V$ is a Harnack curve. Then $f$ is dense i.e., $\supp (f) = \Delta\cap\mathbb{Z}^2$.
\end{The}

\begin{prooof} Theorem 6.3 is a consequence of the second statement of Theorem 6.1, and the fact that the logarithmic map  critical points coincide with the argument map critical points (Passare \cite{P2-06}, see the proof in \cite{N1-08}), and also, the critical values of the logarithmic map coincide with the boundary of the amoeba for Harnack curves (see Mikhalkin \cite{M3-00}).
\end{prooof}

\vspace{0.1cm}

\noindent A description and combinatorial characterization of the indices in $\Ima (\ord_f) \cup \{ type\,\, I\}$ linked to the combinatorial structure of Newton polytope is given in a join work with Petter Johansson \cite{JN-09}.

\vspace{0.2cm}

\section{Examples of complex algebraic plane curves}

\vspace{0.2cm}

In this section, we give some examples of coamoebas of complex algebraic plane curves, and their codual hyperplanes arrangements.

\vspace{0.2cm}

{\bf Example .1-}\,\,
Let $V_{f_{\lambda}}$ be the curve in $(\mathbb{C}^*)^2$ defined by the following polynomial:
$$
f_{\lambda}(z,w) = w^2-\lambda w + 2zw -z^2w +1.
$$
Let $f_{\lambda ,\, 1}(z,w)=  w^2-\lambda w + 2zw -z^2w $, so $V_{f_{\lambda ,\, 1}}$ is just the parabola of example 1.
Let $f_{\lambda ,\, 2}(z,w)=  -\lambda w + 2zw -z^2w +1$, hence $V_{f_{\lambda ,\, 2}}$ is the set of points
$(z,w)\in (\mathbb{C}^*)^2$ such that :
$$
w=\frac{1}{z^2-2z+\lambda}.
$$
This means that $\arg (w_2) = - \arg (w_1)$ modulo $2\pi$.

\begin{figure}[h]
\includegraphics[angle=0,width=0.3\textwidth]{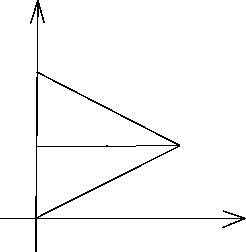}
\caption{Newton Polygon of example 1.}
\label{c}
\end{figure}

Hence the coamoeba of the curve defined by $f_{\lambda}$ is as in figure 14 in the left.

\begin{figure}[h]
\includegraphics[angle=0,width=0.3\textwidth]{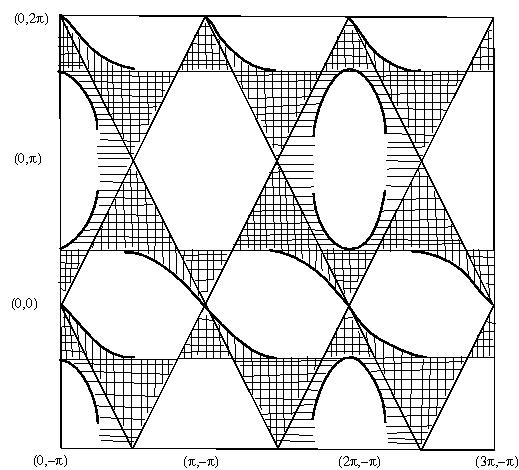}\qquad\qquad\qquad\qquad
\includegraphics[angle=0,width=0.3\textwidth]{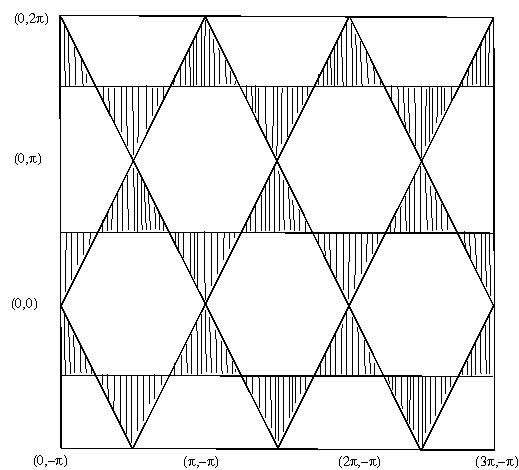}
\caption{Coamoeba of example 1.}
\label{c}
\end{figure}

\vspace{0.2cm}

{\bf Example .2-}\,\,
Let $V_{f_{\lambda}}$ be the curve in $(\mathbb{C}^*)^2$ defined by the following polynomial:
$$
f_{\lambda}(z,w) = zw^2+z^2w+z+w+ \lambda zw.
$$
Let $f_{\lambda ,\, 1}(z,w)= zw^2+z+w+ \lambda zw$. Hence $V_{f_{\lambda ,\, 1}}$ is just
a reparametrization of the parabola of example 1. We can see that $z = -\frac{w}{1+w^2+\lambda w}$.\\
Let $f_{\lambda ,\, 2}(z,w)= zw^2+z^2w+z+ \lambda zw = z(1+zw+w^2+\lambda w) $, hence
$V_{f_{\lambda ,\, 2}}$ is the set of points $(z,w)\in (\mathbb{C}^*)^2$ such that :
$$
z=-\frac{1+w^2+\lambda w}{w}.
$$
This means that $\arg (z_2) = - \arg (z_1)$ modulo $2\pi$, where $z_1$ (resp. $z_2$)
denotes  the first coordinate of a point in $V_{f_{\lambda ,\, 1}}$ (resp. in
$V_{f_{\lambda ,\, 2}}$) . As in example 1, we have the picture in figures 15.

\begin{figure}[h]
\includegraphics[angle=0,width=0.2\textwidth]{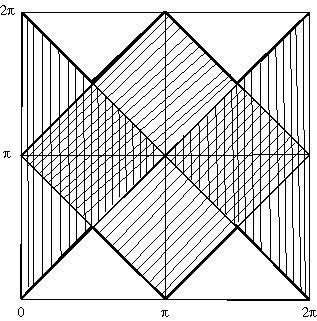}
\includegraphics[angle=0,width=0.4\textwidth]{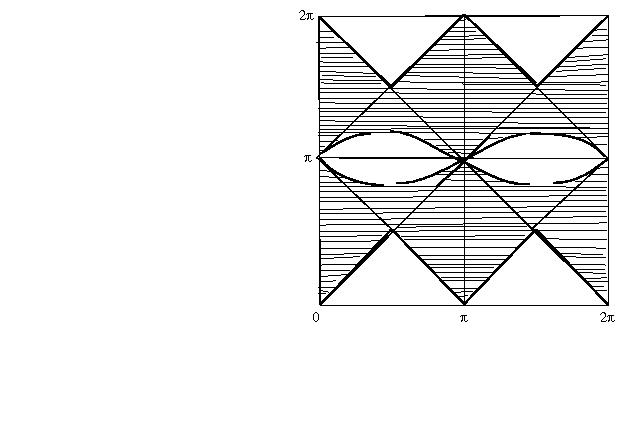}\qquad
\includegraphics[angle=0,width=0.2\textwidth]{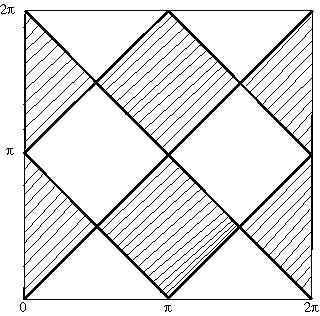}
\caption{Coamoeba of example 2.}
\label{c}
\end{figure}

\vspace{0.2cm}

{\bf Example .3-}\,\,
In this example we illustrate the coamoebas of curves defining by
maximally sparse polynomials with Newton polygon $\Delta$ of
vertices $(0;0),\, (2;1),\, (1;2)$, and $(0;1)$ (see figures 16,
17, and 18).

\begin{figure}[h]
\includegraphics[angle=0,width=0.2\textwidth]{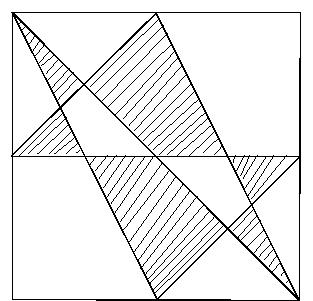}\qquad\qquad
\includegraphics[angle=0,width=0.2\textwidth]{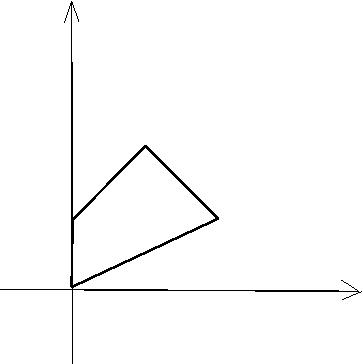}
\caption{Coamoeba of a del Pezzo surface and corresponding triangulation of the Newton polygon dual
to its amoeba.}
\label{c}
\end{figure}

\begin{figure}[h]
\includegraphics[angle=0,width=0.5\textwidth]{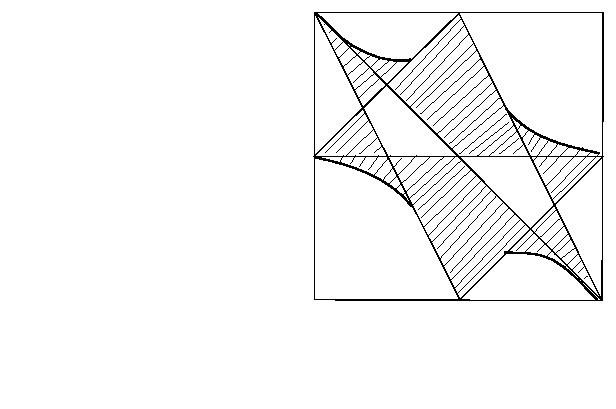}\qquad\qquad
\includegraphics[angle=0,width=0.2\textwidth]{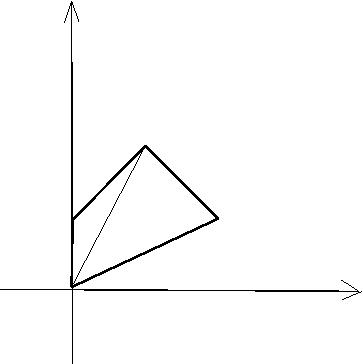}
\caption{Coamoeba of a del Pezzo surface and corresponding subdivision of the Newton polygon dual to
its amoeba.}
\label{c}
\end{figure}

\begin{figure}[h]
\includegraphics[angle=0,width=0.2\textwidth]{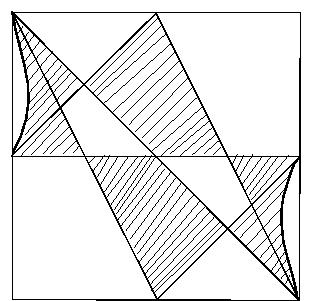}\qquad\qquad
\includegraphics[angle=0,width=0.2\textwidth]{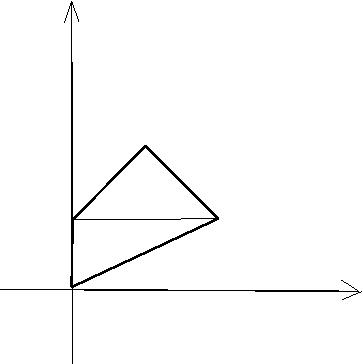}
\caption{Coamoeba of a del  Pezzo surface and corresponding triangulation of the Newton polygon dual
to its amoeba.}
\label{c}
\end{figure}

\begin{figure}[h]
\includegraphics[angle=0,width=0.4\textwidth]{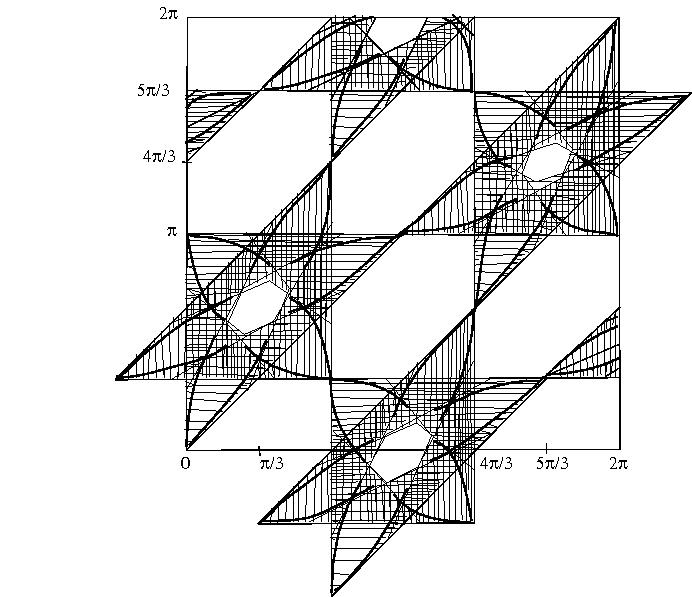}\qquad\qquad
\includegraphics[angle=0,width=0.25\textwidth]{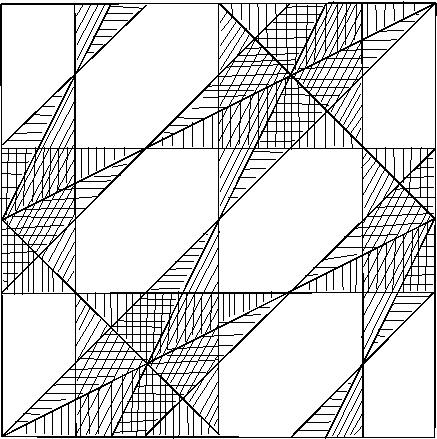}
\caption{Coamoeba of a cubic with solid amoeba on the left and a 
cubic such the spine of its amoeba is of genus one on the right.}
\label{c}
\end{figure}

\vspace{0.2cm}

{\bf Example .4-}\,\,
We give in this section an example of a polygon $\Delta$ such that, there is no real curve defined by a
 polynomial with Newton polygon $\Delta$ and with maximal number of coamoeba complement components, but 
 this maximal number is realized by a complex curve. Let $\Delta$ be the polygon with vertices $(1;0),\,
 (0;1),\, (1;2)$, and $(3;1)$ (see figure 20 for the polygon and its subdivision dual to the spine of the
 amoeba). In this case we prove that no real polynomial can realize the maximal number of coamoeba complement components (the maximal number in the real case is five, and the coamoeba is given in figure 21 on the right for some
 real coefficients ), but the complex curve defined by the  complex polynomial $f(z,w) = 
e^{i\alpha}w+z+zw^2+z^3w$ with $0<\alpha <\pi$, has a coamoeba with maximal number of complement components 
(i.e. six components, see figure 21 on the left; the picture of four fundamenta domains).

\begin{figure}[h]
\includegraphics[angle=0,width=0.2\textwidth]{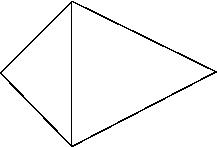}\qquad\qquad\qquad
\includegraphics[angle=0,width=0.2\textwidth]{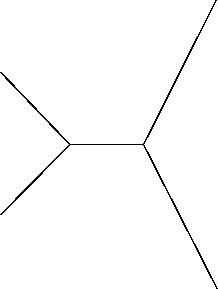}
\caption{The subdivision of the Newton polygon and its dual.}
\label{c}
\end{figure}

\begin{figure}[h]
\includegraphics[angle=0,width=0.4\textwidth]{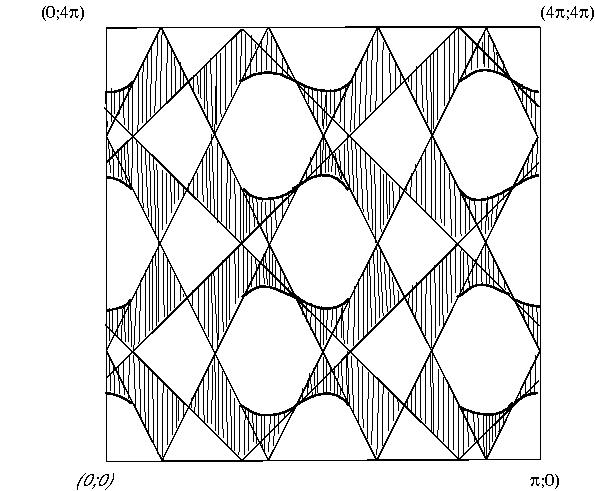}\qquad\qquad
\includegraphics[angle=0,width=0.3\textwidth]{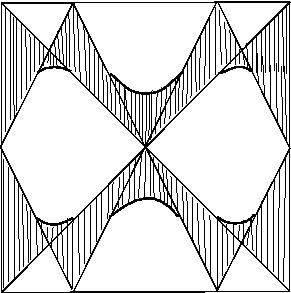}
\caption{Coamoebas of examples 4.}
\label{c}
\end{figure}

\vspace{0.2cm}

{\bf Example .5-} \,\,  We drow in this example the coamoeba of a complex and real parabolas in four fundamental domains (see Figure 22 and 23).

\vspace{0.2cm}

\begin{figure}[h]
\includegraphics[angle=0,width=0.25\textwidth]{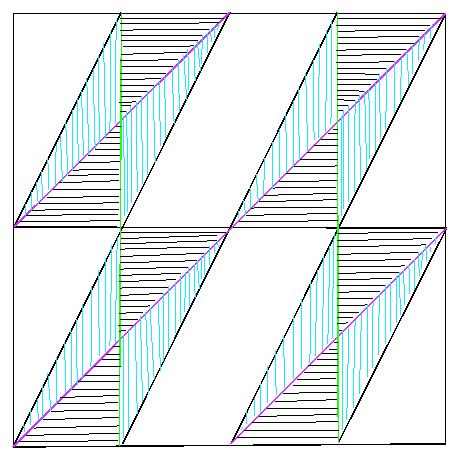}\qquad\qquad
\includegraphics[angle=0,width=0.25\textwidth]{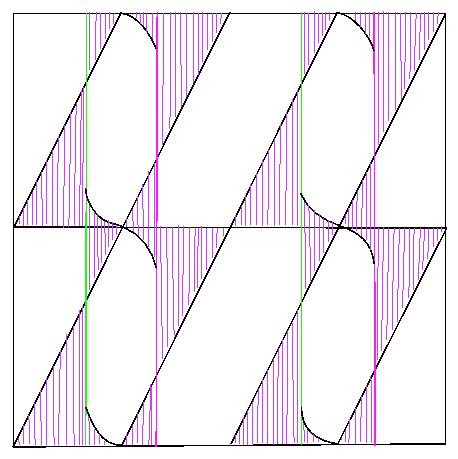}
\caption{The coamoeba of Harnack parabola on the left, and a real parabola with solid amoeba on the right.}
\label{c}
\end{figure}

\begin{figure}[h]
\includegraphics[angle=0,width=0.25\textwidth]{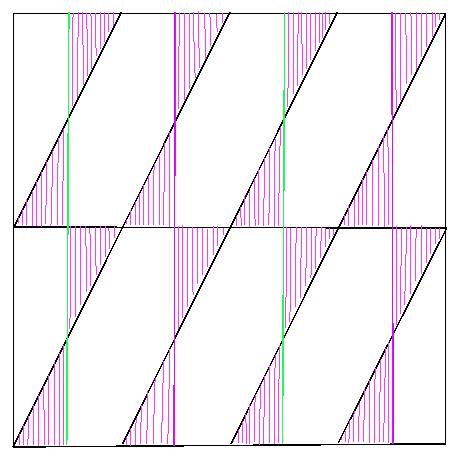}\qquad\qquad
\includegraphics[angle=0,width=0.25\textwidth]{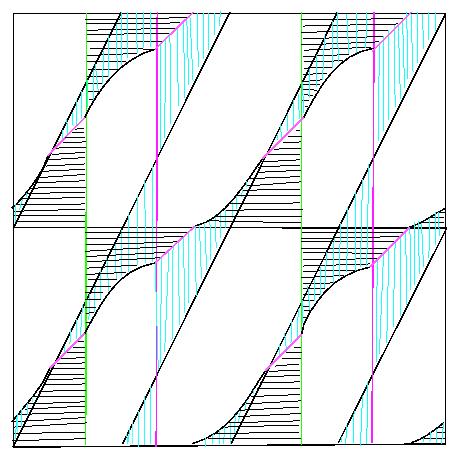}
\caption{Coamoeba of real parabola defined by maximally sparse polynomial on the left, and complex parabola on the right.}
\label{c}
\end{figure}

\end{document}